\documentclass[11pt]{amsart}
\usepackage{geometry}                % See geometry.pdf to learn the layout options. There are lots.
\usepackage{graphicx}
\usepackage{amssymb}
\usepackage{epstopdf}
\usepackage[all]{xy}
\usepackage{stmaryrd}
\newtheorem{thm}{Theorem}[section]
\newtheorem{lem}[thm]{Lemma}
\newtheorem{slem}[thm]{Sublemma}

\newtheorem{cor}[thm]{Corollary}
\newtheorem{fct}[thm]{Fact}

\theoremstyle{definition}
\newtheorem{defn}[thm]{Definition}

\newtheorem{exmp}[thm]{Example}

\newcommand{\revmathfont}[1]{{\textsf{#1}}}
\def\ATRz{{\revmathfont{ATR}$_0$}}
\def\ACAz{\revmathfont{ACA}$_0$}
\def\ACAzp{\revmathfont{ACA}$_0^+$}
\def\RCAz{\revmathfont{RCA}$_0$}
\def\RCAzo{\revmathfont{RCA}$_0^\omega$}
\def\RCAzthr{\revmathfont{RCA}$_0^3$}
\def\ATRz{\revmathfont{ATR}$_0$}
\def\PooCAz{$\Pi^1_1$\revmathfont{-CA}$_0$}
\def\WKLz{\revmathfont{WKL}$_0$}

\def\TSep{$\Sigma^2_1$-\revmathfont{Sep}$^\mathbb{R}$}
\def\DDet{$\Delta_1^\mathbb{R}$-\revmathfont{Det}}
\def\SDet{$\Sigma_1^\mathbb{R}$-\revmathfont{Det}}
\def\CWO{\revmathfont{CWO}$^\mathbb{R}$}
\def\BR{\revmathfont{BR}$_1(\mathbb{R})$}
\def\TR{\revmathfont{TR}$_1(\mathbb{R})$}
\def\SF{\revmathfont{SF}$(\mathbb{R})$}
\def\WO{\revmathfont{WO}$(\mathbb{R})$}
\def\DC{\revmathfont{DC}}
\def\AC{\revmathfont{AC}}
\def\CH{\revmathfont{CH}}

\def\ZF{\revmathfont{ZF}}
\def\ZFC{\revmathfont{ZFC}}

\def\SepClassical{$\Sigma^1_1$-\revmathfont{Sep}}
\def\HYP{\revmathfont{HYP}}
\def\QFAC{\revmathfont{QF-AC}$^{(1, 0)}$}
\def\S{\revmathfont{S}}
\def\K{\revmathfont{K}}

\title{Transfinite Recursion in Higher Reverse Mathematics}
\author{Noah Schweber}
\address{Department of Mathematics, University of California, Berkeley, Berkeley, California 94720 U.S.A. \it E-mail address\rm: schweber@math.berkeley.edu}
\begin{document}

\thanks{The author is grateful to Antonio Montalban and Leo Harrington for numerous helpful comments and conversations. This work will be part of the author's Ph.D. thesis \cite{Schthesis}. The author was partially supported by Antonio Montalban through NSF grant DMS-0901169.}

\begin{abstract} In this paper we investigate the reverse mathematics of higher-order analogues of the theory \ATRz{} within the framework of higher order reverse mathematics developed by Kohlenbach \cite{Koh01}. We define a theory \RCAzthr, a close higher-type analogue of the classical base theory \RCAz, and show that it is essentially a conservative subtheory of Kohlenbach's base theory \RCAzo. Working over \RCAzthr, we study higher-type analogues of statements classically equivalent to \ATRz, including open and clopen determinacy, as well as two choice principles, and prove several equivalences and separations. Our main result is the separation of open and clopen determinacy for reals, using a variant of Steel forcing; in the presentation of this result, we develop a new, more flexible framework for Steel-type forcing.
\end{abstract}

\maketitle
%\section{}
%\subsection{}

\section{Introduction}

The question \[\text{``What role do incomputable sets play in mathematics?"}\] has been a central theme in modern logic for almost as long as modern logic has existed. Six years before Alan Turing formalized the notion of  computability, van der Waerden \cite{vdW30} showed that the splitting set of a field is not uniformly computable from the field; put another way, van der Waerden demonstrated the \it necessity \rm of certain incomputable sets for Galois theory. Other results, especially Turing's solution to the Entscheidungsproblem and the solution by Davis, Matiyasevitch, Putnam, and Robinson of Hilbert's Tenth Problem, established the incomputability of particular sets of natural numbers of interest. In 1975, Friedman \cite{Fri75} initiated the axiomatic study of this question, dubbed ``Reverse Mathematics."

Reverse mathematics requires the choice of both a common language in which to express all analyzed theorems, and a base theory in that language over which all equivalences and non-implications are to be proved. The natural choice of language is that of second-order arithmetic, since it is in this language that computability-theoretic principles are most naturally expressed. The base theory is taken to be \RCAz{}, a precise definition of which is contained in \cite{Sim99}; as a base theory, \RCAz{} is justified by the fact that it captures exactly ``computable" mathematics, in the sense that the $\omega$-models of \RCAz{} are precisely the Turing ideals. One notable feature of reverse mathematics is the existence of the ``Big Five," five subtheories of second-order arithmetic ---  \RCAz, \WKLz, \ACAz, \ATRz, and \PooCAz{} --- each of which is ``robust," in the sense that the same theory results when small changes are made to its exact statement (or to the precise coding mechanisms used), and which correspond to the exact strength, over \RCAz, of the vast majority of theorems studied by reverse mathematics.

However, there is a significant amount of classical mathematics, including parts of measure theory  and most of general topology, which resists any natural coding into the language of second-order arithmetic. This was already recognized by Friedman in \cite{Fri75b}. At the time, the higher-order program failed to draw mathematical attention comparable to that of second-order reverse mathematics.

Recently, however, there has been a return to this question. The framework of finite types --- in which objects of arbitrary finite order, such as sets of sets of reals, are treated directly --- has begun to emerge as a natural setting for a higher reverse mathematics, following Ulrich Kohlenbach's paper on the subject \cite{Koh01}.\footnote{Although it is by no means the only one --- see \cite{Sho13} for an approach via $\alpha$-recursion theory instead, and also \cite{GK13} for a cloesly-related $\alpha$-recursive structure theory. Shore also suggests other approaches which could be interesting, such as via $E$-recursion or the computation theory of Blum-Shub-Smale.} Kohlenbach expands the language of second-order arithmetic to all finite types, and extends the system \RCAz{} to include a version of primitive recursion for arbitrary finite-type functionals. The resulting system, \RCAzo{}, is a proof-theoretically natural conservative extension of \RCAz. (From the point of view of computability theory, however, the choice of base theory may not be so clear; see the discussion at the end of this paper.) 

Work on reverse mathematics in finite types has so far proceeded along one or the other of two general avenues: the analysis of classical theorems about objects not naturally codeable within second-order arithmetic, such as ultrafilters or general topological spaces (\cite{Hun08}, \cite{Kre12}, \cite{Tow11}), or the analysis of higher-type ``uniformizations" of classical theorems of second-order arithmetic (\cite{Koh01}, \cite{SY04}). The present paper instead looks at the higher-type analogues of theorems studied by classical reverse mathematics, focusing in particular on what old patterns hold or fail and what new patterns emerge.\footnote{This is also the approach taken in \cite{Sho13}, there with respect to $\alpha$-recursion rather than finite types.} One natural question in this area is the following: to what extent do the robust subsystems of second-order arithmetic have robust analogues at higher types? 

It is this question which the present paper addresses, focusing on the system \ATRz. In the classical case, much of the robustness of \ATRz{} comes from the fact that being a well-ordering is $\Pi^1_1$-complete. For instance, this is what drives the method of ``pseudohierarchies" by which ill-founded linear orders which appear well-founded, such as those constructed in \cite{Har68}, are used to prove a large number of equivalences at the level of \ATRz; see \cite{Sim99}. Moving up a type, however, changes the situation completely: since we can code an infinite sequence of reals by a single real, the class of well-orderings of subsets of $\mathbb{R}$ is again $\Pi^1_1$, instead of being $\Pi^2_1$ complete. This causes the entire method of pseudohierarchies to break down, and raises doubt that the higher-type analogues of various theorems classically equivalent to \ATRz{} are still equivalent. 

We begin by presenting in section 2 a base theory, \RCAzthr, which is simpler to use than \RCAzo; in section 4, we show that the two theories are essentially equivalent. We then study the complexity over \RCAzthr{} of several higher-type analogues of several principles classically equivalent to \ATRz: comparability of well-orderings, arithmetic bar recursion over a well-founded tree, arithmetic transfinite recursion along an ordinal, clopen determinacy, open determinacy, and $\Sigma^1_1$ separation. In section 2, we show that these principles are almost linearly ordered in terms of strength; the exceptions are due to the need for choice principles when working with higher-type objects, in particular, that the reals are well-orderable and that selection functions for appropriate collections of sets of reals exist. At the bottom of this hierarchy lies the principle asserting the comparability of well-orderings of sets of reals, which we show is remarkably weak at higher types relative to the other principles; at the top, a higher-type version of the separation principle \SepClassical. We also show that the choice principles mentioned above are incomparable.

The main result of this paper, to which section 3 is devoted, concerns the two determinacy principles. In classical reverse mathematics, clopen determinacy fails in \HYP, the model consisting of the hyperarithmetic sets, despite hyperarithmetic clopen games having hyperarithmetic winning strategies, since the method of pseudohierarchies allows us to construct games which are ``hyperarithmetically clopen" but are undetermined in \HYP. This method, as noted above, is no longer valid at higher types, while the complexities of winning strategies for clopen games on reals can still be bounded by a transfinite iteration of an appropriate jump-like operator. This suggests that at higher types, open determinacy becomes strictly stronger than clopen determinacy; using an uncountable version of Steel's tagged tree forcing, we show that this is indeed the case.

The reverse-mathematical results of this paper are summarized in the diagram below:

\bigskip

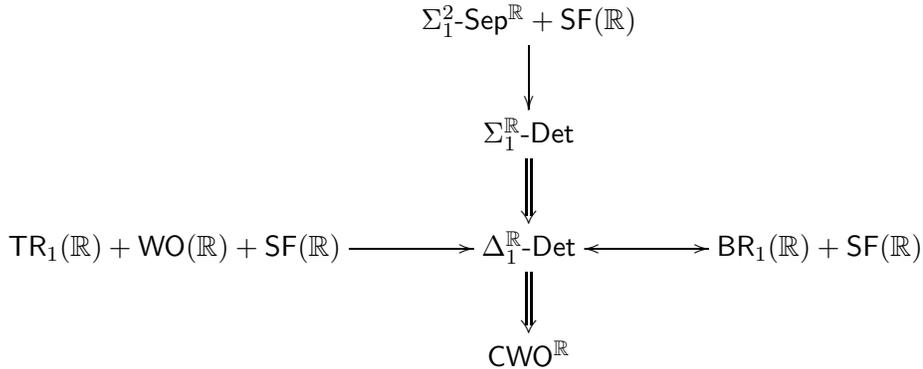
\[\xymatrix{            &  \text{\TSep}+\text{\SF}\ar[d]  &   \\
           &  \text{\SDet} \ar@{=>}[d]    \\
\text{\TR}+\text{\WO}+\text{\SF}\ar@{->}[r]      &  \text{\DDet}  \ar@{=>}[d]  \ar@{<->}[r]  & \text{\BR}+\text{\SF}   \\
               &  \text{\CWO}}\]

\bigskip

Finally, we end by presenting several open questions and directions for further research, raised by both the specific material in this paper and the general program of higher-type reverse mathematics.

\subsection{Background and Conventions}

We refer the reader to \cite{Kun83} for the relevant background in set theory; for descriptive set theory, \cite{Mos09} and \cite{Kec95} are the standard sources. For background on reverse mathematics, see \cite{Sim99}. Finally, for background in finite types, as well as the various computability-theoretic concerns which arise in higher-type settings, see \cite{Lon05}.

\bigskip

There are several notational conventions we adopt for simplicity. Throughout, we use $\mathbb{R}$ to refer to the Baire space, the set of functions from $\omega$ to $\omega$; this is because, during the main result, ordinals will be used as tags, and for this reason a symbol other than ``$\omega^\omega$" is preferable. If $\sigma$ is a nonempty finite string, we write $\sigma^-$ for the immediate $\preccurlyeq$-predecessor of $\sigma$, and if $f$ is an infinite string we write $f^-$ for the string $n\mapsto f(n+1)$.

When writing formulas in many-sorted logic, we use the convention that the first time a variable occurs it is decorated with the appropriate sort symbol; for example, \[\exists x^1\forall y^0(xy=2)\] is the statement ``There is a function from naturals to naturals which is identically 2." (See section 2.1 for a discussion of types.) 

If $\varphi$ is a sentence, then $\llbracket \varphi\rrbracket$ is the truth value of $\varphi$: 1 if $\varphi$ holds, and 0 if $\varphi$ does not. We denote the constant map $k\mapsto 0$ by $\overline{0}$.

If $\Sigma, \Pi\colon A^{<\omega}\rightarrow A$, we write $\Sigma\otimes \Pi$ for the element of $A^\omega$ built by alternately applying $\Sigma$ and $\Pi$: \[\Sigma\otimes \Pi= \langle \Sigma(\langle\rangle), \Pi(\langle \Sigma(\langle\rangle)\rangle), \Sigma(\langle \Sigma(\langle\rangle), \Pi(\langle \Sigma(\langle\rangle)\rangle)\rangle), . . . \rangle.\] We write $(\Sigma\otimes \Pi)_k$ for the length-$k$ initial segment of $\Sigma\otimes \Pi$. A game is said to be a \it win for player X \rm if that player has a winning strategy. A \it quasistrategy \rm for a game played on a set $A$ (so, viewed as a subtree of $A^{<\omega}$) is a multi-valued map from $A^{<\omega}$ to $A$; a quasistrategy is said to be \it winning \rm if each element of $A^\omega$ which is compatible with the quasistrategy is a win for the corresponding player. 

\bigskip

Finally, Theorems \ref{WOvSF} and \ref{Main} rely heavily on the method of set-theoretic forcing. For completeness, we present here a brief summary of this method; for details and proofs, see chapter VII of\cite{Kun83}. 

Given a model $V$ of \ZF{} and a poset $\mathbb{P}\in V$, a \it filter \rm is a subset $F$ of $\mathbb{P}$ which is closed upwards, and such that any two elements of $F$ have a common lower bound in $F$; a set $D\subseteq \mathbb{P}$ is \it dense \rm if every element of $\mathbb{P}$ has a lower bound in $D$.  The \it $\mathbb{P}$-names \rm are defined inductively to be the sets $\{(p_i, \gamma_i): i\in I\}$ of pairs with first coordinate an element of the partial order $\mathbb{P}$, and second coordinate a $\mathbb{P}$-name. If $G$ is a filter meeting every dense subset of $\mathbb{P}$ which is in $V$ --- that is, $G$ is \it $\mathbb{P}$-generic over $V$ \rm --- and $\gamma$ is a $\mathbb{P}$-name, we let $\gamma[G]=\{\theta[G]: \exists p\in G((p, \theta)\in \gamma)\}$ (this is of course a recursive definition). Crucially, the definition of $\gamma[G]$ is made inside $V$, although $G$ will itself will never be in $V$.

We then define the \it generic extension of $V$ by $G$ \rm to be \[V[G]=\{\gamma[G]: \text{$\gamma$ is a $\mathbb{P}$-name in $V$}\}.\] If $V[G]\models\varphi$ whenever $p\in G$, we write $p\Vdash\varphi$; the relation $\Vdash$ is the \it forcing relation \rm given by $\mathbb{P}$. The essential properties of set-theoretic forcing are that the generic extension $V[G]$ is a model of \ZF, and satisfies \AC{} if $V$ does; that the forcing relation is definable in the ground model; and that any statement true in the generic extension is forced by some condition in the generic filter. These are Theorems VII.4.2, VII.3.6(1), and VII.3.6(2) of \cite{Kun83}, respectively.

We will also use the following observation of Fuchs \cite{FHR11} (although we will only need the direction $(1)\implies (2)$) :

\begin{defn} A forcing notion $\mathbb{P}$ is \it countably closed \rm if given any decreasing $\omega$-sequence of conditions \[p_0\ge p_1\ge . . .\] there is some common extension $\in\mathbb{P}$ such that \[q\le p_i, \quad\forall i.\]
\end{defn}

\begin{fct}\label{Count} Let $V\models$\ZF. Then the following are equivalent: \begin{enumerate}
\item $V\models$\DC.
\item Whenever $\mathbb{Q}\in V$ is countably closed in $V$, forcing with $\mathbb{Q}$ adds no new countable sets of ordinals; in particular, forcing with $\mathbb{Q}$ adds no new reals.
\end{enumerate}
\end{fct}

Both forcings considered in this paper are countably closed.

\section{Reverse mathematics beyond type 1}

In this section we first (2.1) review the framework of reverse mathematics in higher types; then (2.2) we define the various higher-type versions of \ATRz{} we will consider in this paper, and prove some basic separations and equivalences. We conclude with Theorem \ref{WOvSF}, whose proof is similar in spirit to, yet much simpler than, that of our main theorem, \ref{Main}.

\subsection{The base theory} 

We begin by making precise the notion of a finite type.\footnote{The one oddity of working with types is that the natural formalization is via many-sorted first-order logic, as opposed to ordinary first-order logic. In many-sorted logic, each element of the model and each variable symbol is labelled by one of a fixed collection of sorts; similarly, function, constant, and relation symbols in the signature must be appropriately labelled with sorts. When there are infinitely many sorts --- as is the case with Kohlenbach's \RCAzo, but not our \RCAzthr{} --- the resulting logic is subtlely different from single-sorted first-order logic; however, these differences shall not be relevant here. For a careful introduction to many-sorted logic, see Chapter VI of \cite{Man96}.} 

\begin{defn}\label{ftypes} The \it finite types \rm are defined as follows: 
\begin{itemize}
\item 0 is a finite type;
\item if $\sigma, \tau$ are finite types, then so is $\sigma\rightarrow \tau$; and
\item only something required to be a finite type by the above rules is a finite type.
\end{itemize}
We denote the set of all finite types by $FT$.

The intended interpretation of finite types is as a hierarchy of functionals, with type 0 representing the ``atomic" objects --- here, natural numbers, or more generally elements of some first-order model of an appropriate theory of arithmetic --- and type $\sigma\rightarrow\tau$ representing the set of maps from the set of objects of type $\sigma$ to the set of objects of type $\tau$.

Within the finite types is the special subclass $ST$ of \it standard finite types\rm,  defined inductively as follows: 0 is a standard type, and if $\sigma$ is a standard type, then so is $\sigma\rightarrow 0$. The standard types are for simplicity identified with natural numbers: $0\rightarrow 0$ is denoted by ``1," $(0\rightarrow 0)\rightarrow 0$ by ``2," etc.
\end{defn}

The appeal of the finite-type framework to reverse mathematics is extremely compelling: the use of finite types lets us talk directly about objects that previously required extensive coding to treat in reverse mathematics, or could not be treated at all. For example, a topological space with cardinality $\le\beth_i$ (where $\beth_0=\aleph_0, \beth_{i+1}=2^{\beth_i}$) can be directly represented as a pair of functionals $(F^i, G^{i+1})$ corresponding to the characteristic functions of the underlying set and the collection of open subsets. Usually, this representation is even natural. In \cite{Koh01}, Kohlenbach developed a base theory for reverse mathematics in all the finite types at once, \RCAzo.

However, working with all finite types at once is cumbersome. First, morally speaking, all finite-type functionals are equivalent to functionals of finite standard type via appropriate pairing functions; second, arbitrarily high types are rarely directly relevant. For that reason, we will use a base theory \RCAzthr, defined below, which only treats functionals of types 0, 1, and 2. In section 4, we prove that the base theory we use in this paper and Kohlenbach's base theory \RCAzo{} are equivalent in a precise sense.

\begin{defn}\label{L3} $L^3$ is the many-sorted first order language, consisting of the following: 

\begin{itemize}
\item Sorts $s_0, s_1, s_2$, with corresponding equality predicates $=_0, =_1, =_2$. We will identify sort $s_i$ with type $i$; recall that the objects of type 0, 1, and 2 are intended to be natural numbers, reals, and maps from reals to naturals, respectively.
\item On the sort $s_0$, the usual signature of arithmetic: two binary functions \[+, \times\colon s_0\times s_0\rightarrow s_0,\] a binary relation \[<\subseteq s_0\times s_0,\] and two constants \[0, 1\in s_0.\]
\item Application operators $\cdot_0, \cdot_1$ with \[\cdot_0\colon s_{1}\times s_0\rightarrow s_0, \quad \cdot_1\colon s_2\times s_1\rightarrow s_0.\] These operators will generally be omitted; e.g., $Fx$ or $F(x)$ instead of $F\cdot_1x$ or $\cdot_1(F, x)$.
\item A binary operation \[*\colon s_2\times s_1\rightarrow s_1\] and a binary operation \[^\smallfrown\colon s_0\times s_1\rightarrow s_1.\]
\end{itemize}

The additional operations $*$ and $^\smallfrown$ allow coding which in Kohlenbach's setting is handled through functionals of non-standard type. Axioms which completely determine $*$ and $^\smallfrown$ are given in Definition \ref{BaseTheory}, below. We will abuse notation slightly and use $^\smallfrown$ to denote both the concatenation of strings, and the specific $L^3$-symbol, as no confusion will arise. Throughout this paper, ``$L^3$-term" will mean ``$L^3$-term with parameters."

Finally, the syntactic classes $\Sigma^0_i$ and $\Pi^0_i$ are defined for $L^3$ as follows:

\begin{itemize}
\item A formula $\varphi$ is in $\Sigma^0_0$ if and only if it has only bounded quantifiers over type 0 objects and no occurrences of $=_1$ or $=_2$. (Note that arbitrary parameters, however, \it are \rm allowed.)

\item A formula $\varphi$ is in $\Pi^0_{i+1}$ if \[\varphi\equiv \forall x^0\theta(x),\] where $\theta\in \Sigma^0_i$.

\item A formula $\varphi$ is in $\Sigma^0_{i+1}$ if \[\varphi\equiv \exists x^0\theta(x),\] where $\theta\in \Pi^0_i$.\end{itemize}

The higher syntactic classes $\Sigma^1_i$, $\Sigma^2_j$, etc. are defined in the analogous way, with lower-type quantifiers being ``for free" as usual.
\end{defn}

The base theory for third-order reverse mathematics which we will use in this paper, \RCAzthr, is then defined as follows: 

\begin{defn}\label{BaseTheory} \RCAzthr{} is the $L^3$-theory consisting of the following axioms: \begin{itemize}
\item $\Sigma^0_1$-induction and the ordered semiring axioms, $P^-$, for the type 0 objects.
\item Extensionality axioms for the type 1 and 2 objects: \[\forall F^1, G^1(\forall x^0(Fx=Gx) \iff F=_1G)\quad\mbox{and}\quad \forall F^2, G^2(\forall x^1(Fx=Gx) \iff F=_2G)\]
\item The $\Delta^0_1$ comprehension\footnote{There are several equivalent formulations of these, including as choice principles; we choose the following presentation, as it is closest to the $\Delta^0_1$ comprehension principle of \RCAzo.} schemes for type 1 and 2 objects:\[\{\forall x^0\exists! y^0\varphi(x, y)\implies \exists f^1\forall x^0(\varphi(x, f(x)))\colon \varphi\in\Sigma^0_1\}\] and \[\{\forall x^1\exists! y^0\varphi(x, y)\implies \exists F^2\forall x^1(\varphi(x, F(x)))\colon \varphi\in\Sigma^0_1\}.\] (The notation ``$\exists!$" is shorthand for ``there exists exactly one.") Recall that $\Sigma^0_1$ formulas may have arbitrary parameters.
\item Finally, the following axioms defining $*$ and $^\smallfrown$: \[\forall k^0, r^1, n^0[(k^\smallfrown r)(n+1)=r(n)\wedge (k^\smallfrown r)(0)=k],\] and
\[\forall F^2, r^1, k^0[(F*r)(k)=F(k^\smallfrown r)].\] The first of these axioms just says that $^\smallfrown$ is the usual concatenation operation; the second is less intuitive, but describes a precise way to turn type-2 functionals into type-$(1\rightarrow 1)$ functionals.
\end{itemize}
\end{defn}

Note that if $(M_0, M_1, M_2; *_0, {}^\smallfrown{}_0), (M_0, M_1, M_2; *_1,  {}^\smallfrown{}_1)\models$\RCAzthr, then in fact \[(M_0, M_1, M_2; *_0,  {}^\smallfrown{}_0)=(M_0, M_1, M_2; *_1,  {}^\smallfrown{}_1);\] that is, models of \RCAzthr{} are determined by their 0-, 1-, and 2-type objects, and it is enough to specify these types to specify the full model. Despite this, the symbols $^\smallfrown$ and $*$ are necessary for \RCAzthr{} since we avoid nonstandard finite types. As evidence of this, the following two facts are easy to prove, yet crucially rely on comprehension over $\Delta^0_1$ formulas involving $^\smallfrown$ and $*$:

\begin{fct} \RCAzthr{} proves each of the following statements:\begin{enumerate}
\item For each type-2 functional $F$, there is a real $r$ such that \[\forall s^1, n^0[\forall k^0(s(k)=n)\implies r(n)=F(s)].\]
\item For each type-2 functional $F$, there is a type-2 functional $G$ such that \[G(\langle a_0, a_1, a_2, . . . , a_n, . . .  \rangle)=F(\langle a_{0}, a_2, a_4, . . . , a_{2n}, . . . \rangle)\]
\end{enumerate}
\end{fct}

\it Proof. \rm For (1), first note that the type-2 comprehension scheme gives us a functional $I$ such that $\forall r^1 [I(r)=r(1)],$ and hence \[\forall r^1, k^0[\forall i^0(I*(k^\smallfrown r)(i)=k)].\] Now our desired real $r$ can be defined by \[r(k)=F(I*(k^\smallfrown \underline{0})),\] which exists by the type-1 comprehension scheme. 

For (2), let $H$ be the type-2 functional defined by the quantifier-free formula $H(r)=r(2r(0)+1)$; then the desired $G$ is defined by the quantifier-free formula \[G(r)=k\iff F(H*r)=k,\] and so again is guaranteed to exist by the type-1 comprehension scheme. $\Box$

\bigskip

It can be shown that neither $(1)$ nor $(2)$ is provable if we restrict the $\Delta^0_1$ comprehension schemes to formulas not involving $*$ and $^\smallfrown$. Essentially, $*$ and $^\smallfrown$ are the price we pay for having axioms which closely resemble those of \RCAz{} and simple models of our base theory. To drive this last point home, we end this section by presenting some natural models of \RCAzthr:

\begin{exmp} Let $\mathcal{I}$ be a Turing ideal. Then the smallest $\omega$-model of \RCAzthr{} containing each real in $\mathcal{I}$ is \[\mathcal{S}_\mathcal{I}=(\omega, \mathcal{I}, \{r\mapsto\varphi^{r\oplus s}_e(0): s\in\mathcal{I}, \varphi_e^{r\oplus s}\text{ total for every $r\in \mathcal{I}$}\}).\] 
\end{exmp}

\begin{exmp}\label{CB} Let \[\mathcal{C}=(\omega, \mathbb{R}, \{F\colon\mathbb{R}\rightarrow\omega\colon F\text{ is continuous}\}), \quad \mathcal{B}=(\omega, \mathbb{R}, \{F\colon\mathbb{R}\rightarrow\omega\colon F\text{ is Borel}\}).\] Both $\mathcal{C}$ and $\mathcal{B}$ are models of \RCAzthr. In both cases, the key step is the observation that \[F_0*(F_1* . . . *(F_n*r))=(\hat{F}_0\circ\hat{F}_1\circ . . . \circ\hat{F}_n)(r)\] (where $\hat{F}\colon r\mapsto\langle F(0^\smallfrown r), F(1^\smallfrown r), . . . \rangle$), and $F$ is continuous/Borel as a function from $\mathbb{R}$ to $\omega$ iff $\hat{F}$ is continuous/Borel as a function from $\mathbb{R}$ to $\mathbb{R}$. This, and the closure of continuous/Borel functions under composition, lets us show that functions of the form \[x\mapsto \llbracket F_0 *(F_1*( . . . *(F_n*x)))(i)=k\rrbracket\] are continuous/Borel so long as their parameters are, which in turn provides the base case for the induction showing $\Delta^0_1$-comprehension for type-2 objects holds in $\mathcal{C}$/$\mathcal{B}$.

We will use both $\mathcal{C}$ and $\mathcal{B}$ in separations in the following subsection (\ref{SFvWO} and \ref{CWO}, respectively).
\end{exmp}

\subsection{Higher-type analogues of \ATRz}

In what follows, we treat higher-type determinacy principles, and towards that end some definitions are necessary. There are several reasonable ways to encode game trees $\subseteq\mathbb{R}^{<\omega}$ as type-2 functionals, and the specific choice of coding is unimportant. When discussing plays, however, things become more complicated. If $\Sigma, \Pi$ are strategies, then the $k$th stage in the play $\Sigma\otimes\Pi$, $(\Sigma\otimes\Pi)_k$ --- or rather, a real coding $(\Sigma\otimes\Pi)_k$ --- is defined as follows. There is a functional $F$, whose existence is guaranteed by the comprehension scheme, such that $F*(k^\smallfrown r)$ is the $k$th ``row" of $r$; specifically, $F$ is defined by \[s\mapsto s(2+\langle s(0), s(1)\rangle).\] We say that a real $r$ \it codes $(\Sigma\otimes\Pi)_k$ \rm if \begin{itemize}
\item $F*(0^\smallfrown r)=\underline{0}$,
\item $\forall 0<2j+1\le k[F*((2j+1)^\smallfrown r)=\Sigma*(F*((2j)^\smallfrown r))]$, and
\item $\forall 0<2j+2\le k[F*((2j+2)^\smallfrown r)=\Pi*(F*((2j+1)^\smallfrown r))]$;
\end{itemize} similarly, we say that $r$ codes the whole play $\Sigma\otimes\Pi$ if $r$ codes $(\Sigma\otimes\Pi)_k$ for all $k$. This definition lets us refer to the play $\Sigma\otimes\Pi$ inside the language of \RCAzthr; and we use, e.g., ``$(\Sigma\otimes\Pi)_k\not\in T$" as shorthand for ``there is a real $r$ coding $(\Sigma\otimes\Pi)_k$, and $r\not\in T$."

There is a subtlety here, which arises due to a particular weakness in the base theory \RCAzo{} (and hence \RCAzthr). (The end of this paper addresses the foundational aspects of this; for now, we simply treat it as it affects us.) \RCAzthr{} is too weak to guarantee the existence of such a real. This is a consequence of Hunter's proof (\cite{Hun08}, Theorem 2.5) that the theory \[\text{\RCAzthr+$\mathcal{E}_1$}:=\text{\RCAzthr+``$\exists E^2\forall x^1(Ex=0\iff \forall k^0(xk=0))$"}\] is conservative over \ACAz: if $\Sigma$ and $\Pi$ are each the operator $E$ described above, then the sentence asserting that $(\Sigma\otimes\Pi)_k$ always exists implies \ACAzp, so that sentence cannot be a consequence of \RCAzthr+$\mathcal{E}_1$.

This can be salvaged in general by altering the base theory; and in fact, since this same subtlety arises in other ways, this is a reasonable course of action --- see the end of Section 5 of this paper. In our case, however, all potential difficulties are handled by the strength of the principles we consider. For example, in the definition of clopen and open determinacy, we use a positive definition of ``winning strategy:" for example, a strategy $\Sigma$ for Open in an open game is winning if for every strategy $\Pi$ for Closed, there is a real coding some fragment $(\Sigma\otimes \Pi)_k$ by which $\Sigma$ has won. This builds into the statements of the theorems we examine all the strength we need to perform the intuitively natural calculations involving stages of games.

The end result is that, although we cannot meaningfully talk about the play of a game $\Sigma\otimes\Pi$ directly within \RCAzthr, the principles we study in this paper happen to have enough power to allow us to do so. As an example of this, it is easy to see that each of the principles defined in Definition \ref{Principles} below imply that at most one player has a winning strategy in an open or clopen game; however, this is not provable in the base theory \RCAzthr{} alone.

\bigskip

Consider the following four theorems, all equivalent to \ATRz{} over \RCAz:
\begin{itemize}
\item \it Comparability of well-orderings\rm: If $X, Y$ are well-orders with domain $\subseteq\mathbb{N}$, then there is an embedding from one into the other.

\item \it Clopen determinacy\rm: Every well-founded subtree of $\omega^{<\omega}$, viewed as a clopen game, is determined.

\item \it Open determinacy\rm: Every subtree of $\omega^{<\omega}$, viewed as an open game, is determined.

\item \it $\Sigma^1_1$ separation\rm: If $\varphi(A)$ is a $\Sigma^1_1$ sentence (possibly with parameters) with a single free set variable, and $X=(X_i)_{i\in\omega}$ is an array of sets such that \[\forall k\in\omega \exists j\in 2(\neg\varphi(X_{\langle k, j\rangle})),\] then there is some set $Y$ such that \[\forall k\in\omega (\neg\varphi(X_{\langle k, Y(k)\rangle})).\]
\end{itemize}

These each have reasonable higher-type analogues, each of which is a theorem of $\text{\ZFC}$:
\begin{defn}\label{Principles} Over \RCAzthr, we define the following principles:
\begin{itemize}

\item The comparability of well-orderings of reals, $\text{\CWO}$: If $X, Y$ are well-orderings with domain $\subseteq\mathbb{R}$, then there is an embedding from one into the other.

\item Clopen determinacy for reals, $\text{\DDet}$: for every tree $T\subseteq\mathbb{R}^{<\omega}$ which is well-founded, viewed as a clopen game, either there is a winning strategy for player I: \[\exists\Sigma\colon \mathbb{R}^{<\omega}\rightarrow\mathbb{R}, \forall \Pi\colon\mathbb{R}^{<\omega}\rightarrow\mathbb{R} [\exists k\in\omega ((\Sigma\otimes\Pi)_{2k+1}\in T\wedge (\Sigma\otimes\Pi)_{2k+2}\not\in T)];\] or there is a winning strategy for player II: \[\exists\Pi\colon \mathbb{R}^{<\omega}\rightarrow\mathbb{R}, \forall \Sigma\colon\mathbb{R}^{<\omega}\rightarrow\mathbb{R}[\exists k\in\omega((\Sigma\otimes\Pi)_{2k}\in T\wedge (\Sigma\otimes\Pi)_{2k+1}\not\in T)].\] 

\item Open determinacy for reals, $\text{\SDet}$: for every tree $T\subseteq\mathbb{R}^{<\omega}$, viewed as an open game, either there is a winning strategy for player I (Open): 

\[\exists\Sigma\colon \mathbb{R}^{<\omega}\rightarrow\mathbb{R}, \forall \Pi\colon\mathbb{R}^{<\omega}\rightarrow\mathbb{R} [\exists k\in\omega ((\Sigma\otimes\Pi)_{k}\not\in T)];\]  or there is a winning strategy for player II (Closed):  \[\exists\Pi\colon \mathbb{R}^{<\omega}\rightarrow\mathbb{R}, \forall \Sigma\colon\mathbb{R}^{<\omega}\rightarrow\mathbb{R}[\forall k\in\omega((\Sigma\otimes\Pi)_k\in T)].\] 

\item The $\Sigma^2_1$-separation principle, \text{\TSep}: If $\varphi(f^2)$ is a $\Sigma^2_1$-formula with a single type-2 free variable, and $X=(X_\eta)_{\eta\in\mathbb{R}}$, $Y=(Y_\eta)_{\eta\in\mathbb{R}}$ are real-indexed collections of type-2 functionals\footnote{A real-indexed set of type-2 functionals $(Z_s)_{s\in\mathbb{R}}$ is coded by the type-2 functional \[\hat{Z}\colon r\mapsto Z_{P_0*r}(P_1*r),\] where $P_0, P_1$ correspond to the left and right projections of a reasonable pairing function $\mathbb{R}^2\cong\mathbb{R}$.} such that \[\neg\exists x^1(\varphi(X_x)\wedge\varphi(Y_x)),\] then there is some type-2 object $F$ such that \[\forall x^1[\varphi(X_x)\implies F(x)=1\quad\mbox{and}\quad \varphi(Y_x)\implies F(x)=0].\]
\end{itemize}

(Note that, strictly speaking, $\text{\TSep}$ is an infinite scheme, as opposed to a single sentence.) It is these principles which we choose to study in this paper. The remainder of this section is devoted to the simpler parts of their analysis; the separation of clopen and open determinacy for reals is the subject of the following section.
\end{defn}

\bigskip

Note that the determinacy principles above are not provable in $\text{\ZF}$ alone, whereas \CWO and $\text{\TSep}$ are, so in order to analyze these principles properly we need some versions of the axiom of choice. The two most important ways choice shows up is in taking the Kleene-Brouwer ordering of a well-founded subtree of $\mathbb{R}^{<\omega}$ and in passing from a quasistrategy for a game on reals to an actual strategy. These correspond, respectively, to the following two principles:

\begin{defn} Letting $\langle\cdot, \cdot\rangle$ be an appropriate pairing function on $\mathbb{R}$, the \it well-ordering \rm and \it selection principles \rm for $\mathbb{R}$ are the following:
\begin{itemize}
\item $\text{\WO}$ is the statement that $\mathbb{R}$ is well-ordered by some functional; that is, there is some type-2 functional $F$ such that the relation \[\{(a, b): F(\langle a, b\rangle)=1\}\] is a well-ordering of $\mathbb{R}$.

\item $\text{\SF}$ is the statement that every $\mathbb{R}$-indexed set of nonempty sets of reals has a selection functional; that is, for every type-2 functional $F$ --- interpreted as the $\mathbb{R}$-indexed set of reals \[\{\{s\in\mathbb{R}: F(\langle r, s\rangle)=1\}: r\in\mathbb{R}\}\] --- there is a type-2 functional $G$ satisfying \[\forall r^1(F(\langle r, G*r\rangle)=1).\]
\end{itemize}
\end{defn}

\bigskip

Now we turn to the implications:

\begin{fct} Over \RCAzthr{}, we have \[\text{\SDet}\implies\text{\DDet}\] and \[\text{\TSep}+\text{\SF}\implies\text{\SDet}.\]
\end{fct}

\it Proof. \rm The first implication is trivial. For the second implication, suppose \TSep{} holds, and let $T\subseteq \mathcal{N}^{<\omega}$ be a tree, viewed as an open game (with Open playing first). We will show $T$ is determined.

Let \[G_\sigma=\begin{cases}
\{\tau: \sigma^\smallfrown\tau\in T\}, & \vert\sigma\vert=2k,\\
\{\tau: \sigma^\smallfrown\tau^-\in T\}, & \vert\sigma\vert=2k+1.\\
\end{cases}\] Note that $G_\sigma$ is a win for Open if and only if $\sigma$ as a node on $T$ is winning for Open; this is the reason for the extra ``padding" move in the case that $\vert\sigma\vert$ is odd. Fixing some reasonable pairing mechanism, let \[H_\sigma=\langle 0, G_\sigma\rangle\quad\mbox{and}\quad J_\sigma=\langle 1, G_\sigma\rangle.\] Let $\varphi$ be the formula \[[\exists \sigma\in T(f=\langle 0, G_\sigma\rangle \wedge \text{ $G$ is a win for Closed})]\vee [\exists \sigma\in T(f=\langle 1, G_\sigma\rangle \wedge \text{ $G$ is a win for Open})].\] We will see that $\varphi$ can be expressed in a $\Sigma^2_1$ way.

First we show how this will prove that $T$ is determined. Note that for a given $\sigma$, we cannot have $\varphi(H_\sigma)$ and $\varphi(J_\sigma)$ both hold, so we can apply \TSep{} (after some appropriate coding). This gives a map \[ev: T\rightarrow \{Closed, Open\}.\] From this map --- and $SF(\mathcal{N})$ --- we can define a winning strategy for $T$: let $\Sigma_{open}$ and $\Sigma_{closed}$ be the strategies for Open and Closed, respectively, which act by selecting (if possible, and playing some fixed real otherwise) an extension to a node on $T$ labelled $Open$ or $Closed$, respectively. It is easy to check that one of these strategies must be winning.

Now to see that $\varphi$ is $\Sigma^2_1$, observe that saying that a game is a win for Open (or Closed, identically) is the same as saying that there is some strategy for Open which defeats all plays. On the face of it this is $\Sigma^2_2$, but each play is just an $\omega$-sequence of reals, and so can be coded as a single real; and now the statement is clearly $\Sigma^2_1$. It follows that $\varphi$, as the disjunction of two $\Sigma^2_1$ statements, is $\Sigma^2_1$. $\Box$

\bigskip

The next implications concern $\text{\DDet}$, and demonstrate its closeness in spirit to the classical second-order system \ATRz:

\begin{defn} We define the following two choiceless versions of $\text{\DDet}$:\begin{itemize}
\item $\text{\BR}$ (bar recursion) is the statement that for any well-founded tree $T\subseteq\mathbb{R}^{<\omega}$, there is a type-2 functional $h$ such that $ran(h)\subseteq\{0, 1\}$ and for all $\sigma\in T$, we have \[h(\sigma)=0\iff \forall a^1(\sigma^\smallfrown\langle a\rangle\in T\implies h(\sigma^\smallfrown\langle a\rangle)=1).\]
(Here we assume some natural coding of elements of $\mathbb{R}^{<\omega}$ by elements of $\mathbb{R}$, so expressions such as $h(\sigma)$ make sense.)
\item $\text{\TR}$ (transfinite recursion) is the scheme of statements asserting that $\Sigma^1_1$ formulas can be iterated along arbitrary well-orderings of sets of reals. More precisely, let $W$ be a well-ordering of a set $dom(W)$ of reals,\footnote{More precisely, we let $W$ be a type-2 functional with range $\subseteq 3$, such that $W(\langle a, b\rangle)=0$ if $a<_Wb$, 1 if $b\le_Wa$, and 2 if $\neg(a, b\in dom(W))$.} and let $\varphi(x^1, Y^2)$ be a $\Sigma^1_1$ formula with the displayed free variables. Then $\text{\TR}$ includes the axiom \[\exists X^2\forall a^1(a\in dom(W)\implies X^{[a]}\equiv ``b^1\mapsto \llbracket \varphi(b, X^{[<_Wa]})\rrbracket"),\] where the functionals $X^{[a]}$ and $X^{[<_Wa]}$ are appropriate codings of the $a$th column and all columns $<_X$-prior to the $a$th column of $W$, respectively.
\end{itemize}
\end{defn}

Unlike \DDet, both \TR{} and \BR{} are provable in \ZF. $\text{\TR}$ is the closest direct analogue of \ATRz, but it also suffers most from the ineffectiveness of the axiom of choice: classically, one passes from transfinite recursion to bar recursion via the Kleene-Brouwer ordering, but taking the Kleene-Brouwer ordering of a subtree of $\mathbb{R}^{<\omega}$ requires a well-ordering of $\mathbb{R}$. \BR{} provides a nice middle ground between determinacy and transfinite recursion, and in fact captures exactly the choiceless part of $\text{\DDet}$.

\begin{lem} Over \RCAzthr{}: 

$(i)$ $\text{\DDet}$ proves $\text{\TR}+\text{\SF}$.

$(ii)$ $\text{\DDet}$ and $\text{\BR}+\text{\SF}$ are equivalent.

$(iii)$ $\text{\TR}+\text{\WO}+\text{\SF}$ proves $\text{\DDet}$.

\end{lem}

\it Proof. \rm For $(i)$, we have $\text{\DDet}\implies \text{\SF}$ easily: given an $\mathbb{R}$-indexed collection of sets of reals, consider the game where player I plays an element of $\mathbb{R}$, and player II then has to play an element of the corresponding set of reals. A winning strategy for this game is exactly the necessary selection function. 

Now let $W$ be a well-ordering of a set $dom(W)$ of reals, and let $\varphi(x^1, Y^2)$ be a $\Sigma^1_1$ formula in the displayed free variables. It is easily checked that \RCAzthr{} proves the following \it countable dependence condition\rm: that for each $a^1, B^2$, if $\varphi(a, B)$ holds then there is a sequence $C=(c_i)_{i\in\omega}$ of reals\footnote{Formally, $C$ is a type-2 functional $F$, and $c_i$ is shorthand for $F*(i^\smallfrown \overline{0})$.} such that \[\forall D^2[\forall i^0(D(c_i)=B(c_i))\implies \varphi(a, D)].\] 

This lets us define the following clopen game: player I plays a real $a_0$, and player II decides whether that real is in the iteration set $X$ required in the definition of $\text{\TR}$. If player II claims that $a_0$ is in $X$, then they have to also supply the values of $X$ on a countable set of reals belonging to columns strictly $W$-below the column of $a_0$; if they claim that $a_0$ is not in $X$, then player I must play such a set. Whichever player did not play the ensuring set, must challenge; this takes the form of asserting that some real $a_1$ is or is not in $X$, where $a_1$ belongs to a $W$-smaller column than $a_0$. A countable set of reals must then be supplied by whichever player claimed $a_1$ is in $X$, and so forth. The first player to be unable to make a legal move, or to play a countable set of reals which do not, in fact, guarantee that some fixed real is in $X$, loses. This is easily checked to be a clopen game, which player I cannot win; and it follows that any winning strategy computes the desired $X$.

For $(ii)$, given a clopen game $G$, use $\text{\BR}$ with the formula \[\varphi(A^2)\equiv\exists a^1(A(a)=0).\] The resulting function $h$ then computes a winning quasistrategy for $G$: if $h(\sigma)=0$, then $\sigma$ is a loss for whoever's turn it is, and one player or the other can win by ensuring that their opponent always plays from nodes marked 0 by $h$. $\text{\SF}$ then lets us pass from this winning quasistrategy to a geniune winning strategy. 

In the other direction, given a well-founded $T$ and a formula $\varphi$, consider the following well-founded game. First, player I chooses some $\sigma\in T$; then, player II responds by playing either ``Safe" or ``Unsafe." The game then continues by playing the clopen game $T_\sigma$, with player II going first if she chose ``Safe" and player II going second if she chose ``Unsafe." Clearly only player II can have a winning strategy, and any winning strategy computes the desired $h$ by setting $h(\sigma)=0$ if the winning strategy for II tells II to play ``Unsafe" if I plays $\sigma$. 

For $(iii)$, the proof is basically the same as $(ii)$, with one slight change: given a clopen game $T$, we use $\text{\WO}$ to form the Kleene-Brouwer ordering of $T$. Then $\text{\TR}$ can be applied, similarly to the use of $\text{\BR}$ above, to get a winning quasistrategy for $T$; finally, we use $\text{\SF}$ to pass to a genuine winning strategy. $\Box$

\bigskip

We now move on to the separations. Among the various forms of \ATRz, there is one particular separation which is straightforward, which we present now. Given the low complexity of wellfoundedness at higher types, it is reasonable to expect that $\text{\CWO}$ is quite weak relative to the higher-type determinacy principles. This is, in fact, true:

\begin{lem}\label{CWO} Let $\mathcal{B}$ be the model generated by the Borel sets, as in Example \ref{CB}. Then \[\mathcal{B}\models\text{\CWO}+\neg\text{\DDet}.\]
\end{lem}

\it Proof. \rm Any uncountable Borel set of reals contains a perfect subset, and there is no Borel well-ordering of $\mathbb{R}$. These facts follow from Borel determinacy (\cite{Kec95}, Theorem 20.5), and together imply that all Borel well-orderings are countable. It then follows that any two Borel well-orderings are comparable by a $\bf\Sigma^0_2$ embedding, so $\mathcal{B}\models \text{\CWO}$.

To show that $B\models\neg\text{\DDet}$, fix some $\bf\Sigma^1_1$ set $X\subseteq \mathbb{R}$ which is not Borel. Let $T\subseteq\omega^{<\omega}$ be a tree such that \[X=\{a\in\mathbb{R}: \exists b\in\mathbb{R} ((\langle a(i), b(i)\rangle)_{i\in\omega}\in [T])\};\] such a tree is guaranteed to exists since $X$ is $\Sigma^1_1$. Now consider the game in which player I plays a real $a$, player II guesses whether $a$ is in $X$, and then the appropriate player (I if II said no, II if II said yes) plays a real $b$ and wins if and only if $b$ is a witness to $a$ being in $X$. This is a clopen game, and viewed as a subtree of $\mathbb{R}^{<\omega}$ it is clearly Borel, and hence in $B$; but if $B$ has a winning strategy for this game, then $B$ must contain $X$, contradicting the assumption that $X$ is not Borel. $\Box$

\bigskip

Note that Borel instances of $\text{\DDet}$ can be constructed whose winning strategies are much more complex than $\Sigma^1_1$; so $\text{\CWO}$ is in fact \it far \rm weaker than $\text{\DDet}$.

\bigskip

Finally, we end this section by establishing the incomparability of the choice principles \WO{} and \SF.

\begin{lem}\label{SFvWO} \RCAzthr+\SF{}$\not\models$\WO.
\end{lem}

\it Proof. \rm Let $C$ be the set of all continuous functions (in some model of \ZFC) from $\mathbb{R}$ to $\omega$; we will see that \[\mathcal{C}=(\omega, \mathbb{R}, C)\models\text{\RCAzthr + \SF + $\neg$\WO.}\] 

Immediately, we have $\mathcal{C}\models\neg$\WO, since there is no continuous well-ordering of the reals; equally immediately, all axioms of \RCAzthr{} except the $\Delta^0_1$-comprehension scheme for type-2 objects hold in $\mathcal{C}$. To show that the comprehension scheme also holds, the key step is showing that any functional defined by a $\Delta^0_1$-formula with continuous functionals as parameters is again continuous; this is an easy yet tedious induction on formula complexity, so we omit it.

Finally, we must show that $\mathcal{C}$ satisfies \SF. To see this, suppose $F$ is an instance of \SF, that is, $F$ is a type-2 functional such that for every real $a$, there is some real $b$ such that $F(\langle a, b\rangle)=1$. Now for $a\in\mathbb{R}$, let $\sigma_a\in\omega^{<\omega}$ be the lexicographically least string such that for all reals $\hat{a}, b$, if $\hat{a}\upharpoonright\vert\sigma_a\vert=a\upharpoonright\vert\sigma_a\vert$ and $\sigma_a\prec b$, then $F(\hat{a}, b)=1$; such a string exists, since $F$ is continuous and is an instance of \SF. More importantly, the map $a\mapsto\sigma_a$ is continuous. From this, it follows that the function \[g\colon\mathbb{R}\rightarrow\omega\colon r\mapsto\begin{cases}
\sigma_{r^-}(r(0)) & \text{ if $r(0)<\vert\sigma_{r^-}\vert$,}\\
0 & \text{otherwise.}
\end{cases}\] is continuous and satisfies $F(a, g(a))=1$ for all reals $a$; so we are done.\footnote{Although this model does have the desired properties, it satisfies \SF{} in a rather unsatisfying way; a more interesting separating model is given by the projective functions, under appropriate large cardinal axioms.} $\Box$

\bigskip

The other separation is more complicated. To the best of our knowledge, there is no natural model of \WO+$\neg$\SF{} as there is of \SF+$\neg$\WO, so we have to build one. Towards this end, we begin with an appropriate model $W$ of \ZF{} in which there is no well-ordering of $\mathbb{R}$, and adjoin a well-ordering of $\mathbb{R}$ by forcing. Of course, this also means that in the generic extension, real-indexed sets of reals have selection functions, so the full model $(\omega, \mathbb{R}, \omega^\mathbb{R})^{W[G]}$ does not separate \SF{} from \WO. Instead, we look at the restricted model \[(\omega, \mathbb{R}^V, \{\nu[G]\colon \nu\in N\})\] for a class $N$ of well-behaved names for type-2 functionals, chosen so that the generic well-ordering of $\mathbb{R}$ winds up in the model, but selection functions for real-indexed nonempty sets of reals do not in general. This is a variation on the basic idea of ``symmetric submodels" which are used to produce models of \ZF{} in which the axiom of choice fails in controlled ways (see \cite{Jec03}, pp. 221-223). The proof of our main result in the following section is also a variation on this basic idea.

\begin{thm}\label{WOvSF} \RCAzthr+\WO$\not\models$\SF.
\end{thm}

\it Proof. \rm We take as the ground model for our forcing argument some \[W\models\text{\ZF+\DC+``The reals are not well-ordered;"}\] the equiconsistency of this theory with \ZF{} itself was proved by Feferman \cite{Fef64}. In $W$, let \[\mathbb{P}=\{p\colon \text{$p$ is a countable partial injective function from $\mathbb{R}$ to $\omega_1$}\}.\] 

First, note that $\mathbb{P}$ is indeed countably closed, so the reals in the generic extension are precisely the reals in $W$. We will use this implicitly in what follows. For $X$ a set, we let $[X]^\omega$ denote the set of countable subsets of $X$. We now define, for $n\in\omega+1$, the \it $n$-countable names \rm inductively as follows:\begin{itemize}
\item A \it 0-code \rm is a pair $c=(c_0, c_1)$, with $c_0\colon \mathbb{R}\rightarrow[\mathbb{R}]^\omega$ and $c_1\colon \mathbb{P}\rightarrow\omega$. If $\nu$ is a name for a map $\mathbb{R}\rightarrow\omega$ and $c$ is a 0-code, we say that $c$ is \it good for $\nu$ \rm if \[\forall p\in\mathbb{P}, a\in\mathbb{R}[c_0(a)\subseteq dom(p)\implies p\Vdash \nu(a)=c_1(p)].\] For $q\in\mathbb{P}$, 0-code $c$ is \it $\nu$-good below $q$ \rm if \[\forall p\le q\in\mathbb{P}, a\in\mathbb{R}[c_0(a)\subseteq dom(p)\implies p\Vdash \nu(a)=c_1(p)].\]
\item Suppose that the set $C_n$ of $n$-codes has already been defined, as well as the notions  ``$\nu$-good" and $\nu$-good below $p$" for $n$-codes. An \it $(n+1)$-code \rm is a pair $c=(c_0, c_1)$ with $c_0\colon\mathbb{R}\rightarrow [\mathbb{R}]^\omega$ and $c_1\colon\mathbb{P}\rightarrow C_n$. If $\nu$ is a name for a map $\mathbb{R}\rightarrow\omega$ and $c$ is an $(n+1)$-code, we say that $c$ is \it $\nu$-good \rm if \[\forall a\in\mathbb{R}, p\in\mathbb{P}[c_0(a)\subseteq dom(p)\implies c_1(p)\text{ is $\nu$-good below $p$}];\] and for $q\in\mathbb{P}$, we say that $c$ is \it $\nu$-good below $q$ \rm if \[\forall a\in\mathbb{R}, p\le q\in\mathbb{P}[c_0(a)\subseteq dom(p)\implies c_1(p)\text{ is $\nu$-good below $p$}].\]
\item A name $\nu$ for a map $\mathbb{R}\rightarrow\omega$ is \it $n$-countable \rm if there is some $n$-code $c$ which is $\nu$-good.
\item $\nu$ is \it $\omega$-countable \rm if $\nu$ is $n$-countable for some $n\in\omega$.
\item Finally, a name $\mu$ for a map $\mathbb{R}\rightarrow\mathbb{R}$ is $n$- or $\omega$-\it countable \rm if the name $\nu$ for the map $r\mapsto \mu(r^-)(r(0))$ is $n$- or $\omega$-countable.
\end{itemize}

The intuition is that the value of an $\omega$-countable name is determined by conditions with large enough domains, mostly regardless of where the elements are sent. This could certainly be pushed past $\omega$, but finite countability is enough for our purposes.

We can now define our target model: Letting $G$ be $\mathbb{P}$-generic over $W$, we set \[M=(\omega, \mathbb{R}^V=\mathbb{R}^{V[G]}, \{\nu[G]: \text{$\nu$ is $\omega$-countable}\}).\] 

Finally, we can finish our proof by showing that $M\models$\RCAzthr+\WO+$\neg$\SF, as follows:

\begin{itemize}
\item $M\models$\WO. This is immediate: the canonical name for the well-ordering \[\prec_G=\{(a, b): G(a)<G(b)\}\] (viewing $G$ as a map $\mathbb{R}\rightarrow\omega_1$) is clearly 0-countable, since to determine whether $G(a)<G(b)$ just depends on $G(a)$ and $G(b)$.

\bigskip

\item $M\not\models$\SF. Our counterexample is $\prec_G$, defined above. Let $\nu$ be $n$-countable. Fix $p\in\mathbb{P}$; we will find a real $a$ and a condition $q\le p$ such that \[q\Vdash \neg (a\prec_G\nu(a)).\] Let $a=\sup (dom(p))+1$, and let $\hat{p}$ be any condition $\le p$ such that $a\in dom(\hat{p})$ and $G(a)-ran(\hat{p})$ is infinite. By induction on $n$, we can ``fill in" the holes in $ran(\hat{p})$ with the reals required to decide $\nu(a)$; that is, we can find $\hat{q}\le\hat{p}$ such that $\sup(ran(\hat{q}))=\hat{q}(a)$, $\hat{q}\Vdash\nu(a)=b$ for some real $b$, and $ran(\hat{q})$ is a \it proper \rm subset of $\hat{q}(a)$. If $b\in dom(\hat{q})$, we take $q=\hat{q}$; if not, we let $q$ be any extension of $\hat{q}$ with $\sup(ran(q))=q(a)$ and $b\in dom(q)$. Either way, the result is a condition, $q$, such that $q\Vdash\nu(a)=b$ but $q(b)<q(a)$, so $\nu$ is not a selection function for $\prec_G$.

\bigskip

\item $M\models$\RCAzthr. All axioms except the $\Delta^0_1$-comprehension scheme for type-2 objects are trivially satisfied, since $M$ is an $\omega$-model containing all the reals. To show that the comprehension scheme holds, note that by a straightforward induction, if $\nu_0$ and $\nu_1$ are $m$- and $n$-countable names for maps $\mathbb{R}\rightarrow\mathbb{R}$ then the name for their composition $\nu_0\circ\nu_1$ is $(m+n)$-countable. From this, it immediately follows that any $\Delta^0_1$ expression $\theta$ with $\omega$-countable parameters defines an $\omega$-countable functional: let $m$ be such that all parameters in $\theta$ are $m$-countable, and let $k$ be the length of $\theta$; then the functional defined by $\theta$ is $mk$-countable.
\end{itemize}

\bigskip

This completes the proof. $\Box$

\section{Separating clopen and open determinacy}

In this section we construct a model $M$ of \RCAzthr{}$+\text{\DDet}+\neg\text{\SDet}$, using a variation of Steel's tagged tree forcing; see \cite{Ste78}, and also \cite{Mon08} and \cite{Nee11}. Throughout this section, we work over a transitive ground model $V$ of $\text{\ZFC}$+\CH. 

The basic picture is as follows. Consider the classic clopen game $G_\alpha$, for $\alpha$ an ordinal, in which players I and II alternately build decreasing sequences of ordinals less than $\alpha$, and the first player whose sequence terminates loses. Clearly player II wins this game, since all she has to do is consistently play slightly larger ordinals than what player I plays.

\bigskip

$G_\alpha\colon \quad$ \begin{tabular}{l|ccccc}
Player I & $\alpha_0$ & & $\alpha_1<\alpha_0$ & & $\cdots$ \\ \hline 
Player II & & $\beta_0$ & & $\beta_1<\beta_0$ & $\cdots$
\end{tabular} $\quad$

\bigskip

Now there is a natural open game, $\mathcal{O}_\alpha$, associated to $G_\alpha$. $\mathcal{O}_\alpha$ has the same rules as $G_\alpha$, except that on player I's turn, she can give up and start over, playing an arbitrary ordinal below $\alpha$. If she does this, then player II gets to play an arbitrary ordinal below $\alpha$ as well. After a restart, play then continues as normal, until player II loses or player I restarts again. Player I (Open) wins if player II's sequence ever reaches zero; player II (Closed) wins otherwise.

\bigskip

$\mathcal{O}_\alpha\colon \quad$ \begin{tabular}{l|ccccc}
Player I (Open) & $\alpha_0$ & & $\alpha_1$ & & $\cdots$ \\ \hline 
Player II (Closed) & & $\beta_0$ & & $\beta_1$ & $\cdots$
\end{tabular} $\quad$ ($\forall i, \alpha_{i+1}<\alpha_i \implies \beta_{i+1}<\beta_i$)

\bigskip

Essentially, $\mathcal{O}_\alpha$ is gotten by ``pasting together" $\omega$-many copies of $G_\alpha$, one after the other, and player II must win \it all \rm of these clopen sub-games in order to win $\mathcal{O}_\alpha$. This is still a win for player II, but in a more complicated fashion. In particular, if player II happened to not be able to directly see the ordinals player I played, but was only able to see the underlying game tree itself, she would need quite a lot of transfinite recursion to be able to figure out what move to play next - seemingly more than she would need to win $G_\alpha$, since there is much more ``noise" in the structure of $\mathcal{O}_\alpha$.

This is essentially the situation we create in the construction below. In a close analogy with Steel forcing we want to force with appropriate partial maps $\subseteq\mathbb{R}^{<\omega}\rightarrow\omega_2\times\omega_2$ to get a ``generic copy" $G$ of $\mathcal{O}_{\omega_2}$ represented on the reals (since we assume \CH{} in $V$); then, by taking the underlying tree of $G$, we can get an open game, $T_G$, classically equivalent to $\mathcal{O}_{\omega_2}$ but where the meaning of the moves --- i.e., the specific ordinals associated to each play --- is hidden. This isn't quite right, since for technical reasons we look instead at partial maps $\subseteq\mathbb{R}^{<\omega}\rightarrow(\omega_2\cup\{\infty\})\times(\omega_2\cup\{\infty\})$, and this corresponds to allowing either player to play $\infty$ as long as she hasn't yet played an ordinal in the same clopen subgame, but the intuition behind the construction remains the same. 

Then, as in Steel forcing, we consider a particular substructure of the full generic extension; however, the approach we take in defining this substructure is essentially the reverse of the usual analysis of Steel forcing. Classically, the desired substructure is defined by first picking out specific elements of the generic extension --- usually paths through a certain tree --- and then closing under hyperarithmetic reducibility; the proof then continues by showing that every element of the resulting model depends only on ``bounded" information about the generic. In our case, we start at the end, and simply consider the part of the generic extension depending on the generic in a ``bounded" way. This is both clearer and more flexible a method than the standard approach; also, higher-type analogue of the hyperarithmetic sets --- the so-called ``hyperanalytic" sets --- is more complicated to work with. See \cite{Mos67} for a definition of this analogue, as well as an account of some early difficulties faced in its study.

By carefully choosing the right notion of boundedness, it turns out we can preserve enough transfinite recursion to get clopen determinacy for reals, $\text{\DDet}$, but are still unable to compute a winnning strategy for $T_G$. One important piece of this picture is that the game $\mathcal{O}_{\omega_2}$, and hence the conditions in the forcing we use, is ``locally clopen" in a precise sense --- this turns out to be necessary to prove the analog of Steel's retagging lemma for this forcing, without which almost nothing can be proved about the resulting model. The other key feature is the countable closure of our forcing notion, a feature necessarily absent in constructions at the second-order level.

\bigskip

The forcing we use in this section is the following:

\begin{defn} Let $\omega_2^*=\omega_2\cup\{\infty\}$, ordered by taking the usual order on $\omega_2$ and setting $\infty>x$ for all $x\in\omega_2^*$ (including $\infty>\infty$). $\mathbb{P}$ is the forcing consisting of all partial maps $p\colon \subseteq \mathbb{R}^{<\omega}\rightarrow \omega_2^*\times\omega_2^*$ satisfying the following conditions, ordered by reverse inclusion: \begin{itemize}
\item $dom(p)$ is a countable subtree of $\mathbb{R}^{<\omega}$ with $p(\langle \rangle)\downarrow=(\infty, \infty)$ (the game starts with player Open moving, and no meaningful tags);
\item $\sigma\in dom(p)\implies [(\vert\sigma\vert=2k+1\wedge p(\sigma^-)_1=p(\sigma)_1) \vee (\vert\sigma\vert=2k\wedge p(\sigma^-)_0=p(\sigma)_0)]$ (player Open is playing $p(\sigma)_0$, Closed is playing $p(\sigma)_1$, and on a given turn exactly one of these values changes); and
\item $\sigma^\smallfrown\langle a, b\rangle\in dom(p)$, $\vert\sigma\vert=2k$, $\infty\not=p(\sigma)_0>p(\sigma^\smallfrown\langle a, b\rangle)_0\implies p(\sigma)_1>p(\sigma^\smallfrown\langle a, b\rangle)_1$ (as long as player Open has not just played an $\infty$, or failed to play less than her previous play, Closed's next play has to be less than her previous play).
\end{itemize}
Note that the way this last condition is phrased allows $p(\sigma)_1$ to be anything when $p(\sigma)_0=\infty$, for $\vert \sigma\vert=2k$, since we have $\infty>\infty$. Also, if  $\vert\sigma\vert=2k$ and $p(\sigma^-)_1=\infty$, then $p(\sigma)_1$ can be anything.

From this point on, we fix a filter $G\subseteq\mathbb{P}$ which is $\mathbb{P}$-generic over $V$.
\end{defn}

The main difference between our forcing $\mathbb{P}$ and Steel forcing is that $\mathbb{P}$ is countably closed (see \ref{Count}). The immediate use of countable closure is that it lets us completely control the type-1 objects in our model; later, we will use countable closure in a more subtle way, to show that no well-orderings of reals of length $\ge\omega_2^V$ are in our model, even though such well-orderings will exist in the full generic extension (Lemma \ref{Bound}).

\bigskip

As with Steel forcing, we have a retagging notion: 

\begin{defn} For $p, q\in\mathbb{P}$ and $\alpha\in\omega_2$, we say that $q$ is an \it $\alpha$-retagging of $p$\rm, and write $p\approx_\alpha q$, if \begin{itemize}
\item $dom(p)=dom(q)$;
\item for $\sigma\in dom(p), i\in 2$ we have \[p(\sigma)_i<\alpha\implies q(\sigma)_i=p(\sigma)_i\] and \[p(\sigma)_i\ge\alpha\implies q(\sigma)_i\ge \alpha.\]
\end{itemize}
\end{defn}

These retagging relations let us define the set of names which depend on the generic in a ``bounded" way:

\begin{defn} Let $\nu$ be a name for a type-2 functional, that is, a map $\mathbb{R}\rightarrow\omega$, and suppose $\alpha\in\omega_2$. Then $\nu$ is \it $\alpha$-stable \rm if for all $a\in\mathbb{R}$, $k\in\omega$, we have \[\forall p, q\in\mathbb{P}[p\approx_\alpha q, p\Vdash \nu(a)=k\implies q\Vdash \nu(a)=k.]\]
\end{defn}

Finally, we can define our desired model:

\begin{defn} Fix $G$ $\mathbb{P}$-generic over $V$. $M$ is defined inductively to be the $L^3$-structure \[M=(\omega, \mathbb{R}, \{\nu[G]: \exists\alpha<\omega_2(\text{$\nu$ is $\alpha$-stable})\}).\]
\end{defn}

The purpose of this section is to prove

\begin{thm}\label{Main} $M\models $\RCAzthr{}$+\text{\DDet}+\neg \text{\SDet}$.
\end{thm}

We begin with two simple properties of the model $M$.

\begin{defn} $T_G$ is the underlying tree of $G$; that is, \[T_G=\{\sigma\in \mathbb{R}^{<\omega}\colon \exists p\in G(\sigma\in dom(p))\}.\]
\end{defn}

\begin{fct}\label{simple} \begin{enumerate}
\item $\mathcal{P}(\omega^\omega)\cap V\subset M_2$.
\item $T_G\in M_2$.
\end{enumerate}
\end{fct}

\it Proof. \rm (1) follows from the fact that canonical names for sets in $V$ do not depend on the poset $\mathbb{P}$, and are hence 0-stable. For (2), the only way to force $\sigma\not\in T_G$ is to have some $p\in G$, $\tau\prec\sigma$ such that $p(\tau)_1=0$, so it follows that the canonical name for $T_G$ is 1-stable. $\Box$

We can now prove the first non-trivial fact about $M$: that it does not satisfy open determinacy for reals. Specifically, we will show that $T_G$, viewed as an open game, is undetermined in $M$. 

The first step is the following:

\begin{lem}\label{game}$V[G]\models T_G$ is a win for Closed.
\end{lem}

\it Proof. \rm By a straightforward density argument, if $G$ is generic, then whenever $\vert\sigma\vert=2k+1$, $p\in G$, and $p(\sigma)_1=\infty$, there is some $q\in G$ and $a\in\mathbb{R}$ such that  $q(\sigma^\smallfrown\langle a\rangle)_1=\infty$. It follows that the strategy \[\Pi(\sigma)=\text{ the $\le_W$-least $a$ such that }\exists p\in G(p(\sigma^\smallfrown\langle a\rangle)_1=\infty)\] is winning for Closed. $\Box$

\bigskip

The indeterminacy of $T_G$ in $M$ then follows from a two-part argument: strategies for Open can be defeated using \ref{game} and the countable closure of $\mathbb{P}$, and stable strategies for Closed can be defeated by pulling the rug out from under her:

\begin{lem}\label{Indet} $M\models\neg\text{\SDet}$.
\end{lem}

\it Proof. \rm Consider the open game corresponding to $T_G$ (in which player I is Open). Recall that $T_G$ is in $M$ and $T_G$ is ``really" a win for player Closed by \ref{simple} and \ref{game}, respectively; we claim that this game is undetermined in $M$.

Suppose $\Sigma$ is a strategy for player Open in $M$. Consider the tree of game-states ``allowed" by $\Sigma$: \[A_\Sigma=\{\sigma\in T_G\colon \exists \Pi(\sigma\prec\Sigma\otimes\Pi)\}.\] Since $T_G$ is actually a win for Closed, the tree $A_\Sigma$ must be ill-founded. Let $f\in V[G]$ be a path through $T_G$. Then $f\in V$, since $\mathbb{P}$ is countably closed and $f$ can be coded by a single real. But then within $V$, we can construct a strategy $\Pi$ which defeats $\Sigma$ by playing along $f$: \[\tau\prec f\implies \Pi(\tau)=f(\vert\tau\vert), \quad\tau\not\prec f\implies\Pi(\tau)=0.\] Since $\Pi$ exists in $V$, $\Pi\in M_2$; so $T_G$ is not a win for Open in $M$.

Now suppose $\Pi$ is a strategy for player Closed in $M$, and suppose (towards a contradiction) that \[p\Vdash \nu\text{ is a winning strategy in }T_G\] where $\nu$ is an $\alpha$-stable name for $\Pi$, $\alpha\in\omega_2$. We can find \begin{itemize} \item $q\le p$, 
\item $a\in\mathbb{R}-\{c\colon \langle c\rangle\in dom(p)\}$, 
\item $b\in\mathbb{R},$ and 
\item $\beta>\alpha$ \end{itemize} such that $\langle a, b\rangle\in dom(q)$, $q(\langle a\rangle)=(\beta, \infty)$, and $q\Vdash\nu(\langle a\rangle)=b$. Now since $q\le p$ and $p$ forces that $\Pi$ wins, we must have $q(\langle a, b\rangle)=(\beta, \gamma)$ with $\gamma>\beta$; so $\gamma>\alpha$. But then we can find a $\hat{q}\approx_\alpha q$ such that $\hat{q}\le p$ and $\hat{q}(\langle a, b\rangle)=(\hat{\beta}, \hat{\gamma})$ for some $\hat{\beta}>\hat{\gamma}$. But then $\hat{q}$ forces that there is some finite play extending $\langle a, b\rangle$ which is a win for Open; and since every possible finite play exists in $M$, this contradicts the assumption that $\nu$ was forced to be a name for a winning strategy. $\Box$

\bigskip

To analyze $M$ further, we require the analogue of Steel's retagging lemma for our forcing:

\begin{lem} \bf (Retagging) \it Suppose $\alpha<\omega_2$ has uncountable cofinality, $p\approx_\alpha q$, $r\le q$, and $\gamma<\alpha$. Then there is some $\hat{r}\le p$ with $\hat{r}\approx_\gamma r$.
\end{lem}

\it Proof. \rm The heart of this proof is the realization that conditions in $\mathbb{P}$, though not well-founded, are ``locally well-founded" in a precise sense. Intuitively, when deciding how to tag a given node of $r'$, we only need to look at a well-founded piece of the domain of $r$; using the ranks of these well-founded pieces as parameters gives us enough ``room" for the natural construction to go through.

Formally, we proceed as follows. Since $\alpha$ has uncountable cofinality, we can find a $\tilde{\gamma}$ such that $\gamma<\tilde{\gamma}<\alpha$ and $\tilde{\gamma}$ is larger than every $r(\sigma)_i$ and $p(\tau)_i$ ($i\in\{0, 1\}, \sigma,\tau\in \mathbb{R}^{<\omega})$ which is less than $\alpha$.

For $\sigma\in dom(r)-dom(p)$, let \[T_\sigma=\{\tau\colon \sigma^\smallfrown\tau\in dom(r)\wedge \forall \rho\prec\tau(\vert\sigma^\smallfrown\rho\vert\text{ odd }\implies \infty\not= r(\sigma^\smallfrown \rho^-)_0>r(\sigma^\smallfrown\rho)_0)\}\] be the set of ways to extend $\sigma$ within $dom(r)$ which according to $r$ don't involve player Open restarting after $\sigma$, and note that for each $\sigma\in dom(r)-dom(p)$ the tree $T_\sigma$ is well-founded. Also, let $N$ be the set of nodes of $dom(r)$ that are new (that is, not in $dom(p)$) but don't follow any new restarts by player Open: \[\{\sigma\in dom(r)-dom(p)\colon \forall \tau\preccurlyeq \sigma( \tau\in dom(r)-dom(p), \vert \tau\vert\text{ odd}\implies r(\tau^-)_0>r(\tau)_0\not=\infty)\}.\] The idea is that we really only need to focus on nodes in $N$: nodes in $dom(p)$ have already had their tags determined, and nodes not in $N\cup dom(p)$ will have no constraints on their tags coming from $p$ at all, since they must follow a restart by Open. In order to define the value of $\hat{r}$ on some node $\sigma$ in $N$, though, we need an upper bound on how large $N$ is above $\sigma$ to keep from running out of ordinals prematurely; this is provided by taking the rank of $T_\sigma$.

Formally, we build the retagged condition as follows. Recalling that $V\models$\ZFC, fix in $V$ a well-ordering of $\mathbb{R}^{<\omega}$, and via that ordering let $rk(S)$ be the rank of $S$ for $S\subseteq \mathbb{R}^{<\omega}$ a well-founded tree. Then we define $\hat{r}$ as follows: 

\[\hat{r}(\sigma)=\begin{cases}
\uparrow, & \text{ if }\sigma\not\in dom(r),\\
p(\sigma), & \text{ if $\sigma\in dom(p)$,}\\
r(\sigma), & \text{ if $\sigma\not\in (N\cup dom(p))$,}\\
(\min\{\tilde{\gamma}+rk(T_\sigma), r(\sigma)_0\}, \hat{r}(\sigma^-)_1), & \text{ if $\sigma\in N$ and $\vert\sigma\vert$ is odd,}  \\
(\hat{r}(\sigma^-)_0, \min\{\tilde{\gamma} + rk(T_\sigma), r(\sigma)_1\}), & \text{ if $\sigma\in N$ and $\vert\sigma\vert$ is even.}\\
\end{cases}\]
It is readily checked that $\hat{r}\in\mathbb{P}$ --- the assumption on $\tilde{\gamma}$ being used here to show that the coordinates of $\hat{r}$ are decreasing when the corresponding coordinates of $r$ drop from $\ge\alpha$ to $<\alpha$ --- and that $\hat{r}\le p$ and $\hat{r}\approx_\gamma r$ (in fact, $\hat{r}\approx_{\tilde{\gamma}} r$). $\Box$

\bigskip

As a straightforward application of the retagging lemma, we can now show that $M$ is a model of \RCAzthr:

\begin{lem}\label{Annoying} $M\models$ \RCAzthr{}.
\end{lem}

$P^-$, the extensionality axioms, the axioms defining $*$ and $^\smallfrown$, and comprehension for reals are all trivially satisfied, the last of these since $M$ contains precisely the reals in $V$ and $V\models$ \ZFC. Only the comprehension scheme for type-2 functionals is nontrivial. We will prove  that \it arithmetic \rm comprehension for type-2 functionals holds in $M$, since this proof is no harder than the proof for $\Delta^0_1$ comprehension. 

Let $\varphi(X^1, y^0)$ be an arithmetic (that is, $\Sigma^0_n$ for some $n\in\omega$; recall Definition \ref{L3}) formula such that for each $a\in \mathbb{R}$ there is precisely one $k\in\omega$ with \[M\models \varphi(a, k).\] Since each natural number is definable, we can assume $\varphi$ has no type-0 parameters. Let $(F_i)_{i<n}$ be the type-2 parameters used in $\varphi$, let $(s_j)_{j<m}$ be the type-1 parameters used in $\varphi$, and let $\nu_i$ be an $\alpha$-ranked name for $F_i$; since each $F_i$ has a ranked name, and there are only finitely many $F_i$, we can find some large enough $\alpha<\omega_2$ so that such names exist. Note that we can work directly with the $s_j$, as opposed to just dealing with their names, since our forcing adds no new reals.

For $a\in\mathbb{R}$, let $\mathcal{C}_a$ be a minimal --- and hence countable --- set of names for reals such that \begin{itemize}
\item $\mathcal{C}_a$ contains a name for $a$ and each $s_j$;
\item whenever a name $\mu$ is in $\mathcal{C}_a$ and $k\in\omega$, $\mathcal{C}_a$ contains a name $\nu$ such that $\Vdash \nu=k^\smallfrown\mu$; and
\item whenever $\mu$ is in $\mathcal{C}_a$ and $i<n$, $\mathcal{C}_a$ contains a name $\mu'$ such that $\Vdash \mu'=\nu_i*\mu$.
\end{itemize}
Although we have not been completely precise in defining the sets $\mathcal{C}_a$, it is clear that the definition above is effective in the sense that a suitable set of sets of names $\{\mathcal{C}_a: a\in\mathbb{R}\}$ exists in the ground model, $V$. Basically, $\mathcal{C}_a$ is the set of reals which can be built from the parameters $s_j$ using the functionals $F_i$ and the operations $*$ and $^\smallfrown$. $\mathcal{C}_a$ is a set of names, as opposed to a set of reals, since the parameters $F_i$ may not exist in the ground model.

The key fact about $\mathcal{C}_a$ is that they determine the truth value of the formula $\varphi$: 

\begin{slem} The truth value of $\varphi(a, k)$ depends only on the values the $F_i$ on the names in $\mathcal{C}_a$. Formally, \[\forall \mu\in\mathcal{C}_a\exists k\in\mathbb{R}[(p\Vdash\mu=k)\wedge (q\Vdash \mu=k)]\implies \forall l\in\omega[(p\Vdash \varphi(a, l))\iff (q\Vdash \varphi(a, l))].\] 
\end{slem}

\it Proof of sublemma. \rm To prove this, first write $\varphi$ as a quantifier-free infinitary formula, with infinite conjunctions and disjunctions in place of universal and existential quantifiers, respectively.

Our claim now follows by induction on the rank of $\varphi$ viewed as a well-founded tree. The induction step is clear; for the base case, the only nontrivial piece is showing that the truth value of a formula of the form \[F_0*(F_1* . . .*(F_m*a))(x)=y\] is determined by the value of the $F_i$ on the (names for) reals in $\mathcal{C}_a$; but this is a straightforward induction on $m$. $\Box$

\bigskip

Now let $\nu$ be a name for the functional defined by $\varphi$. We will show that $\nu$ is $(\alpha+\omega_1)$-ranked.

Let $r\in\mathbb{R}$ and $p, q\in\mathbb{P}$ such that $p\approx_{\alpha+\omega_1}q$ and $p\Vdash \nu(r)=k$. Let \[D_r=\{t\in\mathbb{P}: \forall \mu\in\mathcal{C}_r\exists s\in\mathbb{R}(t\Vdash \mu=s)\}\] be the set of conditions which decide the value of each name in $\mathcal{C}_r$. Since $\mathcal{C}_r$ is countable, and $\mathbb{P}$ is countably closed, the set $D_r$ is dense. Now suppose towards contradiction that $q\not\Vdash \nu(r)=k$. Then since $D_r$ is dense, we can find some $q'\le q$ such that \[q'\in D_r \quad\mbox{and}\quad q'\Vdash \nu(r)=l\] for some natural $l\not=k$. By the retagging lemma, there is some $p'\le p$ such that $p'\approx_\alpha q'$; but since each of the $\nu_i$ are $\alpha$-ranked, we must have \[\forall i<n, t\in\mathbb{R}, \mu\in\mathcal{C}_r[(q'\Vdash \mu=t)\iff (p'\Vdash \mu=t)].\] But since the truth value of $\varphi(r, k)$ depends only on the values of the $\mathcal{C}_r$, this contradicts the fact that $p'\le p$ and $p\Vdash\nu(r)=k$. $\Box$

\bigskip

Showing that $M$ satisfies clopen determinacy for reals, however, requires a more delicate proof. Intuitively, given a stable name for a clopen game, we ought to be able to inductively construct a stable name for a winning (quasi)strategy in that game by just iterating the retagging lemma in the right way. However, since the rank of a stable name is required to be $<\omega_2$, we cannot iterate the retagging lemma $\omega_2$-many times, so we need all clopen games in $M$ to have rank $<\omega_2$. This cannot be derived from the retagging lemma alone; instead, we need to look at particular subposets of $\mathbb{P}$:

\begin{defn} For $\alpha<\omega_2$, $\mathbb{P}_\alpha$ is the subposet of $\mathbb{P}$ defined by \[\mathbb{P}_\alpha=\{p\in\mathbb{P}\colon \forall \sigma\in dom(p), i\in 2(p(\sigma)_i<\alpha \vee p(\sigma)_i=\infty)\}.\]
\end{defn}

Conditions in $\mathbb{P}_\alpha$ will turn out to satisfy a slightly stronger retagging property with respect to $\approx_\alpha$ --- the projecting lemma, below --- than conditions in general, and this will be used to prove that this forcing adds no stable well-orderings of reals longer than any in the ground model. Note that this is false for unstable well-orderings; in particular, forcing with $\mathbb{P}$ collapses $\omega_2$ in the full generic extension.

\begin{defn} For $p\in\mathbb{P}$, $\alpha<\omega_2$, we let the \it $\alpha$-projection of $p$\rm, \[p^\alpha\colon dom(p)\rightarrow (\alpha\cup\{\infty\})\times (\alpha\cup\{\infty\}),\] be the map given by \[\forall \sigma\in dom(p), i\in 2, \quad p^\alpha(\sigma)_i=\begin{cases}
p(\sigma)_i & \text{ if } p(\sigma)_i<\alpha \\
\infty & \text{ otherwise.}
\end{cases}\]
\end{defn}

\begin{lem}\label{Projecting} \bf (Projecting) \it For all $p\in\mathbb{P}$, $\alpha<\omega_2$, we have: 
\begin{enumerate}
\item $p^\alpha\in\mathbb{P}_\alpha;$
\item $p^\alpha\approx_\alpha p;$
\item $p\le q\implies p^\alpha\le q^\alpha;$ 
\item $\vert\mathbb{P}_\alpha\vert^V=\aleph_1;$ and
\item $\mathbb{P}_\alpha$ is countably closed.
\end{enumerate}
\end{lem}

\it Proof. \rm For (1), note that since we set $\infty>\infty$, the map \[x\mapsto\begin{cases}
x & \text{ if $x<\alpha$}\\
\infty & \text{ otherwise}
\end{cases}\] satisfies $x<y\iff \pi(x)<\pi(y)$. So as long as $p$ is in $\mathbb{P}$, the projection $p^\alpha$ will not contain any illegal instances of the second coordinate increasing (which is the only possible obstacle to being a condition), and so will also be in $\mathbb{P}$ - and clearly if $p^\alpha\in\mathbb{P}$, then $p^\alpha\in\mathbb{P}_\alpha$.

(2) and (3) are immediate consequences of (1). Property (3) shows that we can allow $\gamma=\alpha$ in the retagging lemma above if $p$ is assumed to be in $\mathbb{P}_\alpha$, and that we can take $\hat{r}$ to be in $\mathbb{P}_\alpha$ as well in that case. 

For (4), note that elements of $\mathbb{P}_\alpha$ can be coded by countable subsets of $\mathbb{R}\times\omega_1$; the result then follows since $V\models \CH$. 

Finally, for (5), let $(p_i)_{i\in\omega}$ be a sequence of conditions in $\mathbb{P}_\alpha$ with $p_{i+1}\le p_i$. Then since $\mathbb{P}$ is countably closed, we have some $q\in\mathbb{P}$ with $q\le p_i$ for all $i\in\omega$; but then $q^\alpha\in\mathbb{P}_\alpha$ by (1), and since each $p_i\in\mathbb{P}_\alpha$, we have $p_i^\alpha=p_i$ and hence $q^\alpha\le p_i$ by (3). $\Box$

\bigskip

This lemma helps provide us with explicit upper bounds on the lengths of type-2 well-orderings in $M$, via the construction below. We can use this result to provide a bound on the lengths of well-orderings in $M$, which in turn allows the induction necessary for showing clopen determinacy to go through.

\begin{lem}\label{Bound} \bf (Bounding) \it Suppose $\nu$ is a stable name for a well-ordering of $\mathbb{R}$ (that is, $\Vdash\nu$ is a well-ordering of $\mathbb{R}$). Then there is some ordinal $\lambda<\omega_2$ such that \[\Vdash \nu\preccurlyeq\lambda.\] That is, $\omega_2$ is not collapsed in a stable way by forcing with $\mathbb{P}$.
\end{lem}

Suppose $\nu$ is an $\alpha$-stable name for a well-ordering of a set of reals. The proof takes place around the subposet $\mathbb{P}_\alpha$. For a sequence of reals $\overline{a}=\langle a_0, . . . , a_n\rangle$ and a condition $p\in\mathbb{P}$, say that $p$ is \it adequate \rm for $\overline{a}$, and write $Ad(p, \overline{a})$, if $p$ forces that $\overline{a}$ is a descending sequence through $\nu$: \[p\Vdash a_0>_\nu . . . >_\nu a_n.\] Note that since $\nu$ is $\alpha$-stable, $p$ is adequate for $\overline{a}$ if and only if $p^\alpha$ is adequate for $\overline{a}$, by (2) of the previous lemma.

In order to bound the size of $\nu$ in any generic extension, we create in the ground model an approximation to the tree of descending sequences through $\nu$, as follows: \[\mathcal{T}_\nu=\{\langle (p_i, a_i)\rangle_{i<n}\colon p_i\in\mathbb{P}_\alpha\wedge \forall i<j<n(p_j\le p_i\wedge Ad(p_j, \langle a_0, . . . , a_{i-1}\rangle))\}.\] Elements of $\mathcal{T}_\nu$ are potential descending sequences, together with witnesses to their possibility. Now since $\nu$ is a name for a well-ordering, we must have that $\mathcal{T}_\nu$ is well-founded. Otherwise, we would have a sequence of condition/real pairs, $\langle (p_i, a_i)\rangle_{i\in\omega}$,  which build an infinite descending sequence through $\nu$, that is, \[p_{i+1}\le p_i, \quad p_{i+2}\Vdash a_i>_\nu a_{i+1}.\] But then a common strengthening $q\le p_i$, which exists by the countable closure of $\mathbb{P}_\alpha$, would create an infinite descending chain in $\nu$; and this contradicts the assumption that $\Vdash\text{$\nu$ is well-founded}$.

Additionally, $\vert\mathcal{T}_\nu\vert=\aleph_1$, since $\mathcal{T}_\nu\subseteq(\mathbb{P}_\alpha\times\mathbb{R})^{<\omega}$ and $\vert\mathbb{P}_\alpha\vert=\aleph_1$ by Lemma \ref{Projecting}(4). Fixing in $V$ a bijection between $\omega_1$ and $\mathcal{T}_\nu$ we can take the Kleene-Brouwer ordering $\mathcal{L}_\nu$ of $\mathcal{T}_\nu$. Since $\mathcal{T}_\nu$ is well-founded, this is a well-ordering; below, we will show that in fact \[\Vdash \nu\preccurlyeq\mathcal{L}_\nu.\]

Let \[K_\nu^G=\{\langle a_0, . . . , a_n\rangle\colon a_0>_{\nu[G]}. . . >_{\nu[G]}a_n\}\] be the tree of descending sequences through $\nu[G]$ in $V[G]$, and fix a well-ordering $\le_W$ of $\mathbb{P}_\alpha$ in $V$. For $\overline{a}\in K_\nu^G$, we define a condition in $\mathbb{P}_\alpha$ by recursion as follows: \[h(\overline{a})=\text{ the $\le_W$-least $p\in\mathbb{P}_\alpha$ such that $p\le h(\overline{b})$ for all $\overline{b}\prec\overline{a}$ and $Ad(p, \overline{a})$}.\] (Note that by the previous lemma and the fact that $\nu$ is $\alpha$-stable, $h$ is defined for all $\overline{a}\in K_\nu^G$.) An embedding from $K_\nu^G$ into $\mathcal{T}_\nu$ can then be defined: \[e\colon K_\nu^G\rightarrow\mathcal{T}_\nu\colon \langle a_i\rangle_{i<n}\mapsto \langle (h(\langle a_0, . . . , a_i\rangle), a_i)\rangle_{i<n}.\] It follows that $\nu[G]\preccurlyeq\mathcal{L}_\nu$, as desired. $\Box$

\bigskip

Now we are finally ready to prove that $M$ satisfies clopen determinacy. For simplicity, this proof is broken into three pieces. First, we show that the rank of a node in a clopen game can be determined in an $\alpha$-stable way, for appropriately large $\alpha$. Then we define a set which encodes the rank of nodes in a clopen game, as well as which player these nodes are winning for, and show that this set is similarly well-behaved. Finally, we use this to give stable names for winning strategies in clopen games which themselves have stable names --- and this will suffice to show that \DDet{} holds in $M$. Unfortunately, the first two steps in this proof is exceedingly tedious, as we require more and more room to retag conditions, but the intuition is that of a straightforward induction.

\bigskip

Fix in $V$ a well-ordering $\le_W$ of $\mathbb{R}$. Using this well-ordering, we can define the rank $rk(T)$ of a well-founded tree $T\subset\mathbb{R}^{<\omega}$ in the usual way; and for $\sigma\in T$, we let $rk_T(\sigma)=rk(\{\tau: \sigma^\smallfrown\tau\in T\})$. If $\nu$ is a name for a well-founded tree, then $rk(\nu)$ and $rk_\nu(\sigma)$ are the standard names for $rk(\nu[G])$ and $rk_{\nu[G]}(\sigma)$.

\begin{lem} Let $\nu$ be a $\beta$-stable name for a well-founded subtree of $\mathbb{R}^{<\omega}$, $p\in\mathbb{P}$, $\gamma<\omega_2$, and $\sigma\in\mathbb{R}^{<\omega}$ such that \[p\Vdash rk_{\nu}(\sigma)=\gamma,\] and suppose $q\approx_{\beta+\omega_1(\gamma2+2)}p$; then  \[q\Vdash rk_{\nu}(\sigma)=\gamma.\]
\end{lem}

\it Proof. \rm By induction on $\gamma$. For $\gamma=0$, suppose $q$ is a counterexample to the claim; then we can find $r\le q$ and $a\in\mathbb{R}$ such that \[r\Vdash \sigma^\smallfrown\langle a\rangle\in\nu.\] Now by the retagging lemma, we can find some $\hat{r}\le p$ such that $\hat{r}\approx_\beta r$. Since $\nu$ is $\beta$-stable, we have \[r\Vdash \sigma^\smallfrown\langle a\rangle\in\nu,\] which contradicts the assumption on $p$.

Now suppose the lemma holds for all $\gamma<\theta$, and let $p\Vdash rk_{\nu}(\sigma)=\theta$; then \[p\Vdash \forall a\in\mathbb{R}(\sigma^\smallfrown\langle a\rangle\in\nu\implies rk_\nu(\sigma^\smallfrown\langle a\rangle)<\theta).\] Suppose towards a contradiction that \[q\approx_{\beta+\omega_1(\theta2+2)} p \quad\mbox{and}\quad q\not\Vdash rk_{\nu}(\sigma)=\theta;\] then there is some $r\le q$, $a\in\mathbb{R}$ such that \[r\Vdash \sigma^\smallfrown\langle a\rangle\in \nu\wedge rk_\nu(\sigma^\smallfrown\langle a\rangle)\ge\theta.\] By the retagging lemma we get some $\hat{r}\le p$ such that $\hat{r}\approx_{\beta+\omega_1(\theta2+1)}r$, and since $\nu$ is $\beta$-stable we have $\hat{r}\Vdash \sigma^\smallfrown\langle a\rangle\in\nu$. Since $\hat{r}\le p$, and $p\Vdash rk_\nu(\sigma)=\theta$, we must be able to find some $\delta<\theta$ and $s\le \hat{r}$ such that $s\Vdash rk_\nu(\sigma^\smallfrown\langle a\rangle)=\delta$; using the retagging lemma a second time, we can get some $\hat{s}\le r$ such that $\hat{s}\approx_{\beta+\omega_1(\delta2+2)} s$. But then by the induction hypothesis $s\Vdash rk_\nu(\sigma^\smallfrown\langle a\rangle)=\delta$, contradiction the assumption on $r$. $\Box$

\begin{defn} If $T\subset \mathbb{R}^{<\omega}$ is a well-founded tree, thought of as a clopen game, a node $\sigma$ on $T$ is \it safe \rm if the corresponding clopen game \[T^\sigma=\{\tau\colon \sigma^\smallfrown\tau\in T\}\] is a win for player I. For $\nu$ be a $\beta$-stable name for a well-founded subtree of $\mathbb{R}^{<\omega}$ with rank $<\alpha$ for some $\alpha<\omega_2$ (see Lemma \ref{Bound}), let $\Delta_\nu$ be a name for the set which encodes rank and safety of nodes on $\nu$: \[\Delta_\nu[G]:=\{(\sigma, \delta, i): \sigma\in\nu[G]\quad\mbox{and}\quad rk_{\nu[G]}(\sigma)=\delta\quad\mbox{and}\quad i=\llbracket \text{$\sigma$ is safe in $\nu[G]$} \rrbracket\}.\] 
\end{defn}

We will show that $\Delta_\nu$ is well-behaved, in the sense of stability, and use this to give a stable name for a winning strategy for $\nu$.

\begin{lem}\label{DStable} Let $\nu$ be a $\beta$-stable name for a well-founded subtree of $\mathbb{R}^{<\omega}$ of rank $<\alpha$; and for simplicty, let $\kappa=\beta+\omega_1(\alpha2+2)$. If $p\Vdash (\sigma, \delta, i)\in\Delta_\nu$, and $q\approx_{\kappa+ \omega_1(\delta2+2)} p$, then $q\Vdash (\sigma, \delta, i)\in \Delta_\nu$.
\end{lem}

\it Proof. \rm Suppose not. Let $\delta$ be the least ordinal such that for some $\sigma, i$ there are conditions $p, q$ such that \begin{itemize}
\item $q\approx_{\kappa+ \omega_1(\delta2+2)} p$,
\item $p\Vdash (\sigma, \delta, i)\in \Delta_\nu$, and
\item $q\not\Vdash (\sigma, \delta, i)\in \Delta_\nu$.
\end{itemize}

There are two cases. If $\delta=0$, then we must have $i=0$; since $\nu$ is $\beta$-stable, there can be no condition below $q$ which adds a child of $\sigma$ to $\nu$ (since then we can use the retagging lemma to force this below $p$, which already forces that $\sigma$ is terminal in $\nu$), and so $q\Vdash (\sigma, 0, 0)\in \Delta_\nu$.

So suppose $\delta>0$. Since $p\Vdash rk_\nu(\sigma)=\delta$, by the previous lemma we have $q\Vdash rk_\nu(\sigma)=\delta$; so $q$ just disagrees on whether $\sigma$ is safe, which means we must be able to find some $r\le q$ such that \[r\Vdash (\sigma, \delta, 1-i)\in\Delta_\nu.\] By the retagging lemma we can find an $\hat{r}\le p$ such that $\hat{r}\approx_{\kappa+\omega_1(\delta2+1)}r$.

Now the proof breaks into two subcases based on whether $i=0$ or $i=1$.  We treat the first case; the proofs are essentially identical. 

We have $r\approx_{\kappa+\omega_1(\delta2+1)}\hat{r}$ and $r\Vdash (\sigma, \delta, 1)\in \Delta_\nu$. Since $r$ thinks $\sigma$ is safe, $r$ must think there is some immediate successor of $\sigma$ which is unsafe. That is, we can find $s\le r$, $\theta<\delta$, and $a\in\mathbb{R}$ such that $s\Vdash (\sigma^\smallfrown\langle a\rangle, \theta, 0)\in\Delta_\nu$; by retagging again we can find \[\hat{s}\le \hat{r}, \hat{s}\approx_{\kappa+\omega_1(\theta2+2)} s,\] which by our assumption on $\delta$ means that \[\hat{s}\Vdash \sigma^\smallfrown\langle a\rangle\in\nu\text{ and is unsafe}.\] But $\hat{s}\le\hat{r}\le p$ and $p$ believes $\sigma$ is unsafe, which means $p$ believes $\sigma$ has no safe extensions - a contradiction. $\Box$

\bigskip

Finally, we are ready to show that $\nu$ is determined in $M$:

\begin{cor} Let $\nu$ be a $\beta$-stable name for a well-founded subtree of $\mathbb{R}^{<\omega}$, viewed as a clopen game, with $rk(\nu)<\alpha<\omega_2$ for some limit ordinal $\alpha$. Then there is an ($\beta+\omega_1(\alpha4+5)$)-stable name for a (type-2 functional coding a) winning strategy for $\nu$.
\end{cor}

\it Proof. \rm (Note that requiring $\alpha$ to be a limit is a benign hypothesis, as we can always make $\alpha$ larger if necessary; this assumption is just made to simplify some ordinal arithmetic below.) Recall that $\le_W$ is a well-ordering of $\mathbb{R}$ in $V$. Let $\mu$ be a name for the type-2 functional which encodes the strategy picking out the $\le_W$-least winning move at any given stage: \[\mu[G](n^\smallfrown \sigma)=\begin{cases}
a(n) & \text{ if $a$ is the $\le_W$-least real such that $\exists\beta<\alpha[(\sigma^\smallfrown \langle a\rangle, \beta, 0)\in\Delta_\nu$],}\\
0  & \text{ if no such real $s$ exists.}\\
\end{cases}\] For simplicity, we assume that $\Vdash$``no string containing a `0' is on $\nu$," so that there is no ambiguity in this definition. Clearly $\mu$ yields a winning strategy for whichever player wins $\nu$. 

All that remains to show is that $\mu$ is stable. Let $\lambda=(\beta+\omega_1(\alpha4+5))$, fix $\sigma$ and $a$, and let $p\approx_\lambda q$ are conditions in $\mathbb{P}$ such that $p\Vdash \mu(\sigma)=a$. We can find some $r\le q$ and some $b$ such that $q'\Vdash \mu(\sigma)=b$; we'll show that $b=a$, and so we must have had $q\Vdash \mu(\sigma)=a$ already. 

There are two cases:

\bigskip

\it Case 1: $a=0$. \rm Suppose towards a contradiction that $b\not=0$. Since $a=0$, we have $p\Vdash\forall \delta<\alpha, \forall b\in\mathbb{R}[(\sigma^\smallfrown\langle b\rangle, \delta, 0)\not\in\Delta_\nu]$. Let $s\le r$ and $\delta<\alpha$ be such that $s\Vdash (\sigma^\smallfrown\langle b\rangle, \delta, 0)\in \Delta_\nu$; by the retagging lemma, there is $p'\le p$ with $p'\approx_{\beta+\omega_1(\alpha4+4)}s$, which by Lemma \ref{DStable} is impossible.

\bigskip

\it Case 2: $a\not=0$. \rm By identical logic as in the previous case, we must have $b\not=0$; suppose towards contradiction that $b\not=a$. With two applications of the retagging lemma, we can find ordinals $\delta_0, \delta_1<\alpha$ and conditions $p'\le p$, $r'\le r$ such that \begin{itemize}
\item $p'\approx_{\beta+\omega_1(\alpha4+2)} r'$,
\item $p'\Vdash (\sigma^\smallfrown\langle a\rangle, \delta_0, 0)\in\Delta_\nu$, and
\item $r'\Vdash (\sigma^\smallfrown\langle b\rangle, \delta_1, 0)\in\Delta_\nu$.
\end{itemize}
By Lemma \ref{DStable}, we have $r'\Vdash (\sigma^\smallfrown\langle a\rangle, \delta_0, 0)\in\Delta_\nu$ and $p'\Vdash (\sigma^\smallfrown\langle b\rangle, \delta_1, 0)\in\Delta_\nu$ as well. Also note that we have $p'\Vdash\mu(\sigma)=a$, $r'\Vdash\mu(\sigma)=b$, since $p'\le p$ and $r'\le r\le q$. Now since $a\not=b$, either $a<_Wb$ or $b<_Wa$, and so either way we have a contradiction. 

\bigskip

This completes the proof. $\Box$

\bigskip

Since $M_1=\mathbb{R}$, $M$ computes well-foundedness of subtrees of $\mathbb{R}^{<\omega}$ correctly; so by \ref{Bound}, it then follows that every clopen game in $M$ has a winning strategy in $M$. Together with \ref{Indet} and \ref{Annoying}, this completes the proof of Theorem \ref{Main}. $\Box$

\section{\RCAzthr{} versus \RCAzo}

In this section we review Kohlenbach's original base theory \RCAzo, and show that it is equivalent to our base theory \RCAzthr{} in a precise sense.

\begin{defn} Let $L^\omega$ be the many-sorted language consisting of \begin{itemize}
\item a sort $t_\sigma$ for each finite type $\sigma\in FT$,
\item application operators \[\cdot_{\sigma, \rho}\colon t_{\sigma\rightarrow \rho}\times t_\sigma\rightarrow t_\rho\] for all finite types,
\item the signature of arithmetic for the type-0 functionals, and
\item equality predicates $=_\sigma$ for each $\sigma\in FT$.
\end{itemize}
\RCAzo{} is the $L^\omega$-theory consisting of the following axioms: \begin{itemize} 
\item The ordered semiring axioms, $P^-$, for the type-0 objects, and extensionality axioms for all the finite types;
\item the schemata \[\exists \Pi^{\sigma\rightarrow (\tau\rightarrow\tau)}\forall X^\sigma, Y^\tau (\Pi XY=Y)\] and \[\exists \Sigma^{(\sigma\rightarrow (\rho\rightarrow\tau))\rightarrow ((\sigma\rightarrow\rho)\rightarrow (\sigma\rightarrow \tau))} \forall X^{\sigma\rightarrow (\rho\rightarrow \tau)}, Y^{\sigma\rightarrow \rho}, Z^{\sigma}(((\Sigma X) Y) Z=(XZ)(YZ))\] defining the \K- and \S-combinators, respectively;
\item the axiom $\mathcal{R}_0$ asserting the existence of a primitive recursion functional, which for clarity we will always denote $R_0$: \[\exists R_0^{0\rightarrow ((0\rightarrow 1)\rightarrow (0\rightarrow 0))}\forall x^0, g^{0\rightarrow 1}, k^0(R_0(x, g)(0)=x\wedge R_0(x, g)(k+1)=g(R_0(x, g)(k)), k);\] and 
\item the choice scheme \QFAC, which consists --- for each quantifier-free formula $\varphi(X^1, y^0)$ in only the displayed free variables, containing no equality predicate of type $\not=0$ --- of the axiom \[\forall X^1\exists y^0\varphi(X, y)\implies \exists F^2\forall X^1 \varphi(X, F(X)).\]
\end{itemize}
\end{defn}

\bigskip

We will prove that \RCAzo{} is a conservative extension of \RCAzthr.\footnote{Since the language of \RCAzo{} does not include the symbol $*$, it is technically better to say that \RCAzthr{} is a conservative extension of a subtheory of \RCAzo; however, since this will not be an issue, we ignore this point going forward.} 

To begin, we need some basic results about pairing higher-type objects.

\begin{defn} Fix a pairing operator $\langle\cdot, \cdot\rangle$ on natural numbers. In a slight abuse of notation, for reals $x, y$ we let $\langle x, y\rangle$ be the real gotten by pairing $x$ and $y$ pointwise: \[\langle x, y\rangle\colon a\mapsto \langle x(a), y(a)\rangle.\] For a finite sequence $\overline{c}$ of objects which are either all reals or all naturals, let $\langle \overline{c}\rangle$ be the usual coding of $\overline{c}$ by repeated use of the appropriate-type pairing operator $\langle\cdot,\cdot\rangle$, associating to the right.

For $a$ a real, we let $\tilde{a}=a$; for $b\in\omega$, we let \[\tilde{b}=``n\mapsto b."\] For $\overline{c}$ a finite sequence of objects which are each either reals or naturals, let \[\langle \overline{c}\rangle_{_\mathbb{R}}=\langle\tilde{c_0}, \tilde{c_1}, . . . \rangle.\] We write $\pi_i$ for the projection map onto the $i$th coordinate; both \RCAzo{} and \RCAzthr{} prove the existence of the relevant projection functionals. (In the case of \RCAzthr, a real-valued projection $\pi_i(\langle \overline{w} \rangle_{_\mathbb{R}})$ is given by $F_i*(\langle \overline{w}_\mathbb{R} \rangle)$ for a certain functional $F_i$.) Throughout, we use the pairing functions and projection maps in formulas putatively in the language $L^3$ or $L^\omega$, even though those symbols are not in either language, when it is clear that no expressive power is added.
\end{defn}

We begin with the easier result: that \RCAzthr{} is a subtheory of \RCAzo.

\begin{lem} Whenever \[N=(N_\sigma)_{\sigma\in FT}\models \text{\RCAzo},\] we have $(N_0, N_1, N_2)\models$ \RCAzthr{} (with the symbols $^\smallfrown$ and $*$ interpreted in the obvious way).
\end{lem}

\it Proof. \rm The proof that $N$ contains functionals of type $2\rightarrow(1\rightarrow 1)$ and $0\rightarrow (1\rightarrow 1)$ corresponding to $*$ and $^\smallfrown$, respectively, is not hard. It is somewhat tedious, however --- for example, constructing a term corresponding to $^\smallfrown$ requires a definition by cases, and relies on the functional $R_0$ --- and so we omit it.

Since \RCAzo{} proves $P^-$, extensionality, and $\Sigma^0_1$-induction,\footnote{See \cite{Koh01}.} it now suffices to show that $\text{\RCAzo}$ proves the $\Delta_1^0$ comprehension schemata for type 1 and 2 objects; since the former follows in turn from the latter and a bit of coding, we just need to prove $\Delta^0_1$ comprehension for type 2 objects in $\text{\RCAzo}$.

\begin{defn} An $L^3$-formula $\varphi(\overline{x})$ with parameters from $N$ and only type-1 and type-0 variables is \it representable \rm if there is some type-2 functional $F_\varphi\in N$ such that \[N\models F_\varphi(\langle\overline{a}\rangle_{_\mathbb{R}})=1\iff \varphi(\overline{a}).\]
\end{defn}

\begin{slem} All $\Sigma^0_0$ formulas are representable.
\end{slem}

\it Proof of claim. \rm By induction on the number of bounded quantifiers. Note that the representable formulas are closed under negation, so we need only consider one kind of bounded quantifier.

The base case follows immediately from \QFAC: if $\varphi(r)$ has no quantifiers, then \[\psi(r, k)\equiv (k=0\wedge \neg\varphi(r))\vee (k=1\wedge \varphi(r))\] is a quantifier-free formula, and applying \QFAC to $\psi$ yields a representing functional for $\varphi$.

For the induction step, it is enough to show that  \[\varphi(r^1)\equiv\exists x^0<F^2(r)(G^2(\langle x, r\rangle)=1)\] is representable whenever $F, G\in N$ are type-2 parameters. The key tool here is the primitive recursion operator, $R_0$. Using $R_0$, we can define a functional $H$ whose value on a real $r$ is computed by starting with 0, cycling through all naturals less than $F(r)$ and incrementing each time we encounter a solution to $G(\langle -, r\rangle)=1;$ that is, \[H(r)=0\iff \forall x<F(r), G(\langle x, r\rangle)=0.\] Rigorously, we let $H$ be the type-2 functional defined by \[\lambda r^1.R_0(0, \lambda x^0.(\pi_0(x)+G(\langle \pi_1(x), r\rangle)))(F(r)),\] and note that $H$ is clearly in $N$. Then the representing functional we desire is simply \[I\colon r\mapsto \llbracket H(r)>0\rrbracket,\] whose existence in $N$ follows from applying the axiom \QFAC{} to the formula \[\varphi(r, k)\equiv (k=0\wedge H(r)=0)\vee(k=1\wedge H(r)>0).\] This finishes the proof of the sub-lemma. $\Box$

\bigskip

Now we can prove the full $\Delta^0_1$ comprehension scheme in \RCAzo, as follows. Let \[\varphi(X^1, y^0)\equiv \exists z^0\theta(X^1, y^0, z^0)\] be a $\Sigma^0_1$ formula satisfying the hypothesis of the comprehension scheme, with $\theta\in\Sigma^0_0$. By the lemma, let $F\in N$ be the functional such that \[F(\langle A^1, b^0, c^0\rangle_{_\mathbb{R}})=1\iff N\models \theta(A, b, c).\] Now consider the formula \[\psi(X^1, w^0)\equiv w=\langle s, t\rangle\wedge F(\langle X, s, t\rangle_{_{\mathbb{R}}})=1;\] applying \QFAC{} to $\psi$ yields a functional $G\in N$ of type 2, and $\varphi$ is represented by the functional \[X^1\mapsto \pi_0(G(X)),\] which is clearly in $N$. $\Box$
 
\bigskip

The other half of the equivalence is a conservativity result:

\begin{thm}\label{Conserve} $\text{\RCAzo}$ is conservative over $\text{\RCAzthr}$, in the following sense: given any model $(M_0, M_1, M_2; ^\smallfrown, *)$ of \RCAzthr, there is a model $N=(N_\sigma)_{\sigma\in FS}$ of \RCAzo{} with the same first-, second-, and third-order parts and corresponding application operators: \[M_0=N_0, M_1=N_1, M_2=M_2, \quad \cdot_{0}^M=\cdot_{0, 0}^N, \cdot_{1}^M=\cdot_{1, 0}^N.\]
\end{thm}

\it Proof. \rm Let $M=(M_0, M_1, M_2)\models$ \RCAzthr. Define the set of \it $\lambda$-terms over $M$ \rm as follows:

\begin{defn}\label{Lambda} Fix $M\models$\RCAzthr. The set of \it $\lambda$-terms over $M$ \rm is defined inductively as follows:\begin{itemize}
\item If $\mathfrak{t}$ is an $L^3$-term of type $\sigma$, then $\mathfrak{t}$ is a $\lambda$-term of type $\sigma$.
\item If $\mathfrak{t}, \mathfrak{s}$ are $\lambda$-terms of types $\sigma\rightarrow\tau$ and $\sigma$ respectively, then $\mathfrak{((t)(s))}$ is a $\lambda$-term of type $\tau$.
\item If $\mathfrak{t}$ is a $\lambda$-term of type $\sigma$, $x$ is a variable of type $\rho$, and the expression ``$\lambda x$" does not occur in $\mathfrak{t}$, then $\lambda x.(\mathfrak{t})$ is a $\lambda$-term of type $\rho\rightarrow\sigma$.
\end{itemize}
Free and bound variables are defined as usual. A $\lambda$-term $\lambda x^\sigma. \theta(x)$ is intended to denote the map $a\mapsto \theta(a)$, and so the set of (appropriate equivalence classes of) $\lambda$-terms is meant to define a type-structure. For simplicity, we will refer to $\lambda$-terms over $M$ simply as ``$\lambda$-terms."

We let $T$ be the set of all $\lambda$-terms. We say that $\mathfrak{t}\in T$ is \it closed \rm if $\mathfrak{t}$ has no free variables, that is, if each variable appearing in $\mathfrak{t}$ is within the scope of a $\lambda$.

Now, let $\equiv_{ext}$ be the smallest equivalence relation on $T$ satisfying the following: \begin{itemize}
\item If $a, b$ are $L^3$-terms and $M\models a=b$, then $a\equiv_{ext}b$.
\item For  all $\lambda$-terms $\mathfrak{s}$ of type $\sigma$ in which the variable $x^\sigma$ does not appear and all $\lambda$-terms $\mathfrak{t}$, we have \[((\lambda x^\sigma. \mathfrak{t})(\mathfrak{s}))\equiv\mathfrak{t}[\mathfrak{s}/x].\] (Recall that $\mathcal{e}[a/b]$ denotes the expression gotten by replacing each occurence of $a$ by $b$ in $\mathcal{e}$.)
\item If $\mathfrak{t},\mathfrak{s}$ are $\lambda$-terms of type $(\sigma\rightarrow\tau)$ such that \[\mathfrak{t}(\mathfrak{a})\equiv_{ext}\mathfrak{s}(\mathfrak{a})\] for all $\mathfrak{a}$ of type $\sigma$, then $\mathfrak{t}\equiv_{ext}\mathfrak{s}$.
\item If $\mathfrak{t}, \mathfrak{s}$ are $\lambda$-terms of type 1 and $\mathfrak{t}(a)\equiv_{ext}\mathfrak{s}(a)$ for all $a\in M_0$, then $\mathfrak{t}\equiv_{ext}\mathfrak{s}$.
\item If $\mathfrak{t}, \mathfrak{s}$ are $\lambda$-terms of type 2 and $\mathfrak{t}(a)\equiv_{ext}\mathfrak{s}(a)$ for all $a\in M_1$, then $\mathfrak{t}\equiv_{ext}\mathfrak{s}$.
\end{itemize}
\end{defn}

Our desired model, $N$, will be built out of these $\equiv_{ext}$-classes. The first step towards an analysis of $\equiv_{ext}$-classes of $\lambda$-terms is the following classical result, here stated in a form slightly weaker than usual but more directly useful for our purposes:

\begin{defn} A $\lambda$-term $\mathfrak{t}$ is in \it normal form \rm if it contains no subterm of the form $((\lambda x. \mathfrak{s})(\mathfrak{u}))$
\end{defn}

\begin{thm}\label{Normal} \bf (Normal Form Theorem) \it For each $\lambda$-term $\mathfrak{t}$, there is a $\lambda$-term $\mathfrak{s}$ in normal form such that \[\mathfrak{t}\equiv_{ext}\mathfrak{s}.\]
\end{thm}

See, e.g., section 4.3 of \cite{GTL89} for a proof. The value of the normal form theorem is that it allows us to focus our attention on only nicely-behaved $\lambda$-terms. The relevant nice behavior is captured in the following lemma:

\begin{lem}\label{NiceForm} Let $\mathfrak{t}$ be a $\lambda$-term in normal form. Then: \begin{enumerate}
\item Every subterm of $\mathfrak{t}$ is in normal form.
\item If $\mathfrak{t}$ has standard type (that is, type 0, 1, 2, etc.), then every subterm of $\mathfrak{t}$ also has standard type.
\item If $\mathfrak{t}$ is of type 0 or 1, then all bound variables in $\mathfrak{t}$ are of type 0.
\item If $\mathfrak{t}$ is of type 2, then $\mathfrak{t}$ contains at most one bound variable of type 1, and all other bound variables are of type 0.
\item If $\mathfrak{t}$ is of type 0, 1, or 2, then every subterm of $\mathfrak{t}$ is of type 0, 1, or 2.
\end{enumerate}
\end{lem}

\it Proof. \rm $(1)$ is immediate. For $(2)$, suppose that $\mathfrak{t}$ has standard type, and suppose $\mathfrak{t}$ has subterms of non-standard type. Note that, by induction, every $\lambda$-term of non-standard type has form either $\mathfrak{s}_0(\mathfrak{s}_1)$ or $\lambda x^\sigma. \mathfrak{s}_0$. Now let $\mathfrak{s}$ be the minimal leftmost $\lambda$-subterm of $\mathfrak{t}$ which has nonstandard type; that is, let $\mathfrak{s}$ be the unique subterm of $\mathfrak{t}$ such that $(i)$ no subterm of $\mathfrak{t}$ containing any characters to the left of $\mathfrak{s}$ is of nonstandard type, and $(ii)$ $\mathfrak{s}$ is the shortest subterm of $\mathfrak{t}$ with property $(i)$. Then $\mathfrak{s}$ clearly cannot be of the form $\mathfrak{s}_0(\mathfrak{s}_1)$, since then $\mathfrak{s}_0$ would also need to be of nonstandard type and would then contradict the minimality of $\mathfrak{s}$ among leftmost nonstandard-type $\lambda$-terms; so $\mathfrak{s}$ is of the form $\lambda x^\sigma.\mathfrak{u}$. But since there is no $\lambda$-term of nonstandard type occuring to the left of $\mathfrak{s}$, and only $\lambda$-terms of nonstandard type can be applied to $\mathfrak{s}$ on the left, we must have that $\mathfrak{s}$ is bound by an application on the right. This immediately contradicts the normality of $\mathfrak{t}$, by (1).

For $(3)$, suppose $\mathfrak{t}$ contains a bound variable $y^\sigma$ of type $\sigma\not=0$. The subterm in which $y$ appears bound, $\lambda y. \mathfrak{s}$, then has type $\sigma\rightarrow\tau$ for some $\tau\in FT$. As in the proof of (2), let $\mathfrak{u}$ be the minimal leftmost subterm of $\mathfrak{t}$ which contains $\lambda y.\mathfrak{s}$ and has type $\not=0$; since $\mathfrak{t}$ is in normal form, $\mathfrak{u}$ must be of the form \[\lambda x_0^{\sigma_0}.\lambda x_1^{\sigma_1}. . . \lambda x_n^{\sigma_n}.\lambda y^\sigma.\mathfrak{s}.\] Now since $\mathfrak{t}$ has type 0 or 1, $\mathfrak{u}$ must be on the left or right side of an application. Since $\mathfrak{t}$ is in normal form, $\mathfrak{u}$ must be on the right side of an application; that is, $\mathfrak{t}$ contains a subterm of the form $((\mathfrak{v})(\mathfrak{u}))$. But $\mathfrak{v}$ cannot be a parameter from $M$, since the type of $\mathfrak{u}$ has height at least 2, so $\mathfrak{v}$ is of the form $\lambda z.\mathfrak{w}$; so $((\mathfrak{v})(\mathfrak{u}))$ is not in normal form, contradicting (1).

$(4)$ follows similarly to $(3)$. Since $\mathfrak{t}$ has type 2, $\mathfrak{t}$ may have the form $\lambda y^1. \mathfrak{s}$, in which case it contains at least one bound variable of type 1; but then $\mathfrak{s}$ has type 0, and so by $(3)$ $y$ is the only  bound variable of type $\not=0$ occuring in $\mathfrak{t}$. 

$(5)$ follows the pattern of the previous parts. Towards contradiction, consider the minimal leftmost subterm $\mathfrak{s}$ of $\mathfrak{t}$ of type $n>2$; then no functional can apply to $\mathfrak{s}$ on the left, and binding $\mathfrak{s}$ with a $\lambda$ would result in a subterm of nonstandard type, contradicting $(2)$, so $\mathfrak{s}$ must be on the left of an application. But by minimality $\mathfrak{s}$ has the form $\lambda x^1. \mathfrak{u}$, so this contradicts the normality of $\mathfrak{t}$. $\Box$

\bigskip

\begin{defn} Let $T^*_i$ be the set of all closed $\lambda$-terms of type $i$, for $i\in\{0, 1, 2\}$. For $i\in 3$, let $e_i\colon M_i\rightarrow T^*_i$ be the map \[e_i\colon a\mapsto [a]_{ext}.\]
\end{defn}

Our immediate goal is to show that the $e_i$ are bijections. The remainder of the conservativity proof will then follow easily. Injectivity is straightforward; to show surjectivity, we use the following construction:

\begin{defn} Fix a type-1 variable $y^1$. Let $\mathcal{T}$ be the set of ordered pairs $(\mathfrak{t}, S)$, where $\mathfrak{t}$ is a $\lambda$-term in normal form of type 0 or 1 containing no free variables of nonzero type besides possibly $y$, and $S=(x_0, . . . , x_m)$ is a list of type-0 variables including all those occurring in $\mathfrak{t}$.

For $(\mathfrak{t}, (x_i)_{i<n})\in\mathcal{T}$, say that $F\in M_2$ \it codes \rm $(\mathfrak{t}, (x_i)_{i<n})$ if \begin{itemize}
\item either $\mathfrak{t}$ is of type 0, and for all $a_0, . . . , a_{n-1}\in M_0, b\in M_1$, we have \[F(a_0{}^\smallfrown. . . {}^\smallfrown a_{n-1}{}^\smallfrown b)\equiv_{ext} \mathfrak{t}(a_0/x_0, . . . , a_{n-1}/x_{n-1}, b/y);\]
\item or $\mathfrak{t}$ is of type 1, and for all $a_0, . . . , a_{n-1}\in M_0, b\in M_1$, we have \[F*(a_0{}^\smallfrown. . . {}^\smallfrown a_{n-1}{}^\smallfrown b)\equiv_{ext} \mathfrak{t}(a_0/x_0, . . . , a_{n-1}/x_{n-1}, b/y).\]
\end{itemize}
\end{defn}

The key lemma is the following:

\begin{lem}\label{Associates} Every $(\mathfrak{t}, (x_i)_{i<n})\in \mathcal{T}$ is coded by some $F\in M_2$.
\end{lem}

\it Proof. \rm We will prove the lemma by induction on the complexity of $\mathfrak{t}$; note that by Lemma \ref{NiceForm}, if $(\mathfrak{t}, (x_i)_{i<n})\in \mathcal{T}$, then $(\mathfrak{t}, (x_i)_{i<n})\in\mathcal{T}$ for all subterms $\mathfrak{s}$ of $\mathfrak{t}$, so an induction is possible. For this induction, we make the following abbreviations. For $\pi$ a permutation of $n$ for $n\in\omega$, we let $R_\pi$ be the type-2 functional satisfying 
\[\forall r\in M_1, i\in M_0\colon (R_\pi*r)(i)=\begin{cases} r(\pi(i)) & \text{ if $i<n$,}\\
r(i) & \text{ if $n\le i$.}\\
\end{cases}\] For $n\in\omega$, we let $P_n$ be the type-2 functional satisfying \[\forall r\in M_1, i\in M_0\colon (P_n*r)(i)=r(n+i).\] The existence of such functionals in $M_2$ is an easy consequence of the type-2 comprehension scheme. Finally, recall (\ref{L3}) the definition of the language $L^3$, as well as our convention that ``$L^3$-term" means ``$L^3$-term with parameters." The induction then proceeds as follows:

\bigskip

Fix $(\mathfrak{t}, (x_i)_{i<n})\in\mathcal{T}$, and suppose that for all subterms of $\mathfrak{t}$, and all appropriate lists of variables, the result holds.

For $\mathfrak{t}$ an $L^3$-term, the comprehension scheme for type-2 functionals gives us the desired $F$ immediately. If $\mathfrak{t}$ has type 0, apply comprehension to the formula \[\Phi(u^1, v^0)\equiv v=\mathfrak{t}[u(0)/x_0, . . . , u(n-1)/x_{n-1}, (P_n*u)/y],\] and if $\mathfrak{t}$ has type 1, apply comprehension to the formula \[\Psi(u^1, v^0)\equiv v=\mathfrak{t}[u(1)/x_0, . . . , u(n)/x_{n-1}, (P_{n+1}*u)/y](u(0)).\] Clearly the so-defined $F$ codes $(\mathfrak{t}, (x_i)_{i<n})$.

If $\mathfrak{t}$ is of the form $\mathfrak{s}_0+\mathfrak{s}_1$, by induction let $F_0, F_1$ represent $(\mathfrak{s}_0, (x_i)_{i<n})$ and $(\mathfrak{s}_1, (x_i)_{i<n})$ respectively. Then \[F\colon u^1\mapsto F_0(u)+F_1(u),\] whose existence is again guaranteed by comprehension, clearly codes $(\mathfrak{t}, (x_i)_{i<n})$. (Multiplication and successor are handled identically.)

If $\mathfrak{t}=\lambda z^0. \mathfrak{s}$, note that $\mathfrak{t}$ necessarily has type 1. By induction let $G$ represent $(\mathfrak{s}, (z, x_0, . . . , x_{n-1})).$ Then $G$ also codes $(\mathfrak{t}, (x_i)_{i<n})$.

If $\mathfrak{t}=\mathfrak{s_0}(\mathfrak{s_1})$, then by \ref{NiceForm}(5) $\mathfrak{s_0}$ has type either 1 or 2. If $\mathfrak{s}_0$ has type 1 (and so $\mathfrak{s}_1$ has type 0),  let $F_0$ and $F_1$ represent $(\mathfrak{s}_0, (x_i)_{i<n})$ and $(\mathfrak{s}_1, (x_i)_{i<n})$ respectively. Then \[F\colon u^1\mapsto (F_0*u)(F_1(u))\] codes $(\mathfrak{t}, (x_i)_{i<n})$, and is guaranteed to exist by comprehension.

Finally, suppose $\mathfrak{t}=\mathfrak{s_0}(\mathfrak{s_1})$ and $\mathfrak{s}_0$ has type 2. Note that $\mathfrak{s}_0$ cannot be of the form $\lambda y^1. \mathfrak{u}$, since $\mathfrak{t}$ is in normal form. Similarly, if $\mathfrak{s}_0$ were of the form $\mathfrak{u}(\mathfrak{v})$, then $\mathfrak{u}$ would have to have non-standard type, and this would contradict \ref{NiceForm}(2). This leaves as the only possibility that $\mathfrak{s}_0$ is a single type-2 parameter, that is, $\mathfrak{t}$ is of the form \[F(\mathfrak{s}_1)\] for some parameter $F\in M_2$. By induction, let $G$ code $\mathfrak{s}_1$; then the functional $H$ defined by \[H\colon r\mapsto F(G*r)\] is in $M$ by comprehension, and codes $\mathfrak{t}$.

Since these are all the cases which can arise, this completes the induction. $\Box$

\bigskip

As an immediate corollary, we get the surjectivity of $e_2$: 

\begin{cor} For each closed $\lambda$-term $\mathfrak{t}$ of type 2, there is some $F\in M_2$ such that $F\equiv_{ext}\mathfrak{t}$.
\end{cor}

\it Proof. \rm By \ref{Normal}, we can assume $\mathfrak{t}$ is in normal form; then either $\mathfrak{t}$ is an $L^3$-term, in which case we are done, or $\mathfrak{t}$ has the form \[\mathfrak{t}=\lambda y^1. \mathfrak{s}.\] Applying Lemma \ref{Associates} to $\mathfrak{s}$, we get an $F\in M_2$ such that for all $b\in M_1$, $F(b)=\mathfrak{s}[b/y]$; and so $\mathfrak{t}\equiv_{ext} F$. $\Box$

\bigskip

This then passes to $e_1$ and $e_0$:

\begin{cor} For each closed $\lambda$-term $\mathfrak{t}$ of type 0 (type 1), there is an $a\in M_0$ ($b\in M_1$) such that $\mathfrak{t}\equiv_{ext}a$ ($\mathfrak{t}\equiv_{ext}b$).
\end{cor}

\it Proof. \rm For the type 1 case, let $\lambda x^0. \mathfrak{s}$ be a $\lambda$-term of type 1; now consider the type-2 term $\mathfrak{t}'=\lambda y^1. \mathfrak{s}[y(0)/x]$; by Corollary 4.14, $\mathfrak{t}'\equiv_{ext}F$ for some $F\in M_2$. But then consider the real $f\colon k^0\mapsto F(k^\smallfrown \overline{0})$, which is in $M_1$ by comprehension; clearly $f\equiv_{ext}\mathfrak{t}$ since $F\equiv_{ext}\mathfrak{t}'$, so we are done.

The type 0 case follows similarly. Let $\mathfrak{t}$ be a $\lambda$-term of type 0, and consider the type-2 term $\mathfrak{t}'=\lambda y^1. \mathfrak{t}$. Taking $F\in M_2$ such that $F\equiv_{ext}\mathfrak{t}'$, we must have $F(\overline{0})\equiv_{ext}\mathfrak{t}$; but $F(\overline{0})\in M_0$. $\Box$

\bigskip

Now consider the following $L^{\omega}$-structure $N_M$, defined as follows:\begin{itemize}
\item $N_\sigma=T^*_\sigma$.
\item Application in $N$ is defined as \[\cdot_{\sigma, \rho}\colon ([\mathfrak{t}]_{ext}, [\mathfrak{s}]_{ext})\mapsto[((\mathfrak{t})(\mathfrak{s}))]_{ext}.\]
\item The arithmetic functions are transferred from $M_0$ to $N_0$ in the obvious way; for example, \[+^N\colon ([\mathfrak{a}]_{ext}, [\mathfrak{b}]_{ext})\mapsto [\mathfrak{a}+\mathfrak{b}]_{ext},\] etc.
\item The relation $<^N$ is defined by \[<^N=\{([a]_{ext}, [b]_{ext}): a, b\in M_0, M\models a<b\}.\]
\end{itemize}

\begin{lem} $N_M\models$\RCAzo.
\end{lem}

\it Proof. \rm Extensionality and $P^-$ for the type-0 functionals hold trivially. The $\Pi$ and $\Sigma$ combinators are easily expressed as $\lambda$-terms, and so the corresponding axioms are satisfied.

The primitive recursion axiom, $\mathcal{R}_0$, takes a bit more work to express as a $\lambda$-term. By $\Sigma^0_1$ induction in $M$, for any real $r^1\in M_1$ and any naturals $a^0, b_0\in M_0$, there is a (natural number coding a) primitive recursive derivation for \[R_0(a, r)(b)=k\] for a unique $k\in M_0$. This lets us apply the $\Delta^0_1$ comprehension scheme for type-2 functionals, and so it follows that there is an $F\in M_2$ such that $F(a^0{}^\smallfrown b^0{}^\smallfrown r^1)=k^0$ if and only if there is a code for a primitive recursive derivation of $R_0(r, a)(b)=k$. Now consider the $\lambda$-term \[\mathfrak{t}:= \lambda x^0.\lambda r^1. \lambda y^0. (F(x^\smallfrown y^\smallfrown r));\] by definition of $F$, $[\mathfrak{t}]_{ext}$ clearly witnesses the axiom $\mathcal{R}_0$ in $N$.

Finally, the choice principle \QFAC requires a bit of work: an appropriate quantifier-free formula $\Phi$ may contain high-type parameters, in which case the comprehension scheme in $M$ does not directly apply. Instead, we have to essentially lower the types of the parameters, using the coding provided by Lemma \ref{Associates}.

Let $\Phi(y^1, x^0)$ be a quantifier-free formula in the displayed free variables, containing no occurrences of $=_\sigma$ for $\sigma\not=0$. Then since $\Phi$ contains no higher-type equality predicates, every maximal term in $\Phi$ must have type 0. Let $\mathfrak{t}_0, . . . , \mathfrak{t}_m$ be a list of these maximal terms, so that $\Phi$ is a Boolean combination of formulas of the form $\mathfrak{t}_i=\mathfrak{t}_j$ or $\mathfrak{t}_i<\mathfrak{t}_j$ for $i, j\le m$. 

Note that each $\mathfrak{t}_i$ has free variables from among $\{y^1, x^0\}$. Thus, by Lemma \ref{Associates}, we can find functionals $F_0, . . . , F_m\in M_2$ such that $F_i(a^\smallfrown b)\equiv_{ext}\mathfrak{t}_i[a/x, b/y]$ for every $a\in M_0, b\in M_1$. Let $\hat{\Phi}$ be the quantifier-free formula gotten from $\Phi$ by replacing each $\mathfrak{t}_i$ by $F_i(x^\smallfrown y)$ for each $i\le m$; then applying comprehension to $\hat{\Phi}$ yields a type-2 functional $G\in M_2$, and it is easily checked that \[\forall b\in N_1(\Phi(b, [G]_{ext}(b))).\]

This completes the proof. $\Box$

\bigskip

Theorem \ref{Conserve} is thus proved. $\Box$

\section{Conclusion}

In this paper, we have sought to understand how the passage to higher types affects mathematical constructions related to the system \ATRz; given both the sheer number of such constructions, and the relative youth of higher-order reverse mathematics, this remains necessarily incomplete. We close by mentioning four particular directions for further research we find most immediately compelling.

\subsection{}

Despite the analysis provided by this paper, there are still basic questions remaining unaddressed. It is unclear whether \SDet{} implies \TSep, or whether \DDet{} implies \WO. We suspect that neither implication holds, but separations at this level are unclear: for example, it is open even whether \SDet{} implies the $\Pi^2_2$-comprehension principle for type-2 functionals, although the answer is almost certainly that it does not.

For that matter, in this paper we have focused entirely on the strengths of third-order theorems relative to other third-order theorems; their strength relative to second-order principles has been completely unexplored. For instance, it is entirely possible, albeit unlikely, that \DDet{} and \SDet{}  have the same second-order consequences.

\subsection{}

There is also a lingering question regarding our main theorem: can clopen and open determinacy for reals be separated in a simpler way? Specifically, there is a reasonable ``smallest natural" model $N=(\omega, \mathbb{R}^L, N_3)$ of \RCAzthr+\DDet, with $N_3$ the set of all type-2 functionals in $L_\alpha$ where $\alpha$ is the least ordinal such that this model satisfies $\text{\DDet}$. It is then natural to ask whether $N\models$\SDet. Either answer to this question would be interesting, as well as answering this question for other natural models of \RCAzthr+\DDet.

\subsection{}

One interesting aspect of the shift to higher types we have not touched on at all is the extra structure available in higher-type versions of classical theorems.  Given a $\Pi^1_2$ principle \[\varphi\equiv\forall X^1\exists Y^1\theta(X, Y),\] we can take its higher-type (so prima facie $\Pi^2_2$) analogue \[\varphi^*\equiv\forall F^2\exists G^2\theta^*(F, G).\] Now, individual reals are topologically uninteresting, but passing to a higher type changes the situation considerably. Specifically, we can consider \it topologically restricted \rm versions of $\varphi^*$: given a pointclass $\Gamma$, let \[\varphi^*[\Gamma]\equiv \forall F^2\in\Gamma\exists G^2\theta(F, G).\] The relevant example is restricted forms of determinacy: the principles $\text{\DDet}[\Gamma]$ (resp., $\text{\SDet}[\Gamma]$) assert determinacy for clopen (resp., open) games whose underlying trees when viewed as sets of reals are in $\Gamma$. In particular, the system $\text{\SDet}[Open]$ is extremely weak, at least by the standards of higher-type determinacy theorems: it is equivalent over \RCAzthr{} to the classical system \ATRz. 

The techniques used in the proof of Theorem \ref{Main} are topologically badly behaved. In particular, they tell us nothing about the restricted versions $\text{\DDet}[\Gamma]$ and $\text{\SDet}[\Gamma]$. With some work the argument of this paper might extend to showing that \DDet$[\Gamma]\not\vdash\text{\SDet}[\Gamma]$ over \RCAzthr{}, for reasonably large pointclasses $\Gamma$, but not immediately; and certainly a detailed understanding of which restricted forms of open determinacy for reals are implied by which restricted forms of clopen determinacy will require substantially new ideas. This finer structure seems to allow a rich connection between classical descriptive set theory and higher reverse mathematics, and is worth investigating.

\subsection{}

Finally, there is a serious foundational question regarding the base theory for higher-order reverse mathematics. The language of higher types is a natural framework for reverse mathematics, as explained at the beginning of section 2.1; however, the specific base theory \RCAzo{} is not entirely justified from a computability-theoretic point of view. While proof-theoretically natural, it does not necessarily capture ``computable higher-type mathematics." The most glaring exmaple of this concerns the Turing jump operator. In the theory \RCAzo{}, the existence of a functional corresponding to the jump operator \[\mathfrak{J}^{1\rightarrow 1}\colon f\mapsto f'\] is conservative over \ACAz{} (\cite{Hun08}, Theorem 2.5). However, intuitively we can compute the $\omega$th jump (and much more) of a given real by iterating $\mathfrak{J}$; thus, given a model $M$ of \RCAzo, there may be algorithms using only parameters from $M$ and effective operations which compute reals not in $M$. From a computability theoretic point of view, then, \RCAzo{} may be an unsatisfactorily weak base theory.

Of course, this discussion hinges on what, precisely, ``computability" means for higher types. A convincing approach is given in \cite{Kle59}, justified by arguments by Kleene and others (see especially \cite{Gan67}) similar in spirit to Turing's original informal argument. It is thus desirable --- at least for higher-type reverse mathematics motivated by computability theory, as opposed to proof theory --- to have a base theory corresponding to full Kleene recursion.\footnote{It should be noted that the separation of $\text{\DDet}$ and $\text{\SDet}$ in this paper does not suffer from the choice of base theory. This is because --- by an induction on indices --- Kleene computability from a type-2 object satisfies the following \it countable use condition\rm: if $F$ is a given type-2 object, and $\varphi_e^F$ is a type-2 object computed from $F$, then for each real $r$ there is a countable set of reals $C_r$ such that \[\forall G^2(G\upharpoonright C_r=F\upharpoonright C_r\implies \varphi_e^F(r)=\varphi_e^G(r)).\] It follows by essentially the same argument as \ref{Annoying} that $M$ is closed under Kleene reduction. In fact, all separations in this paper --- with the possible exception of \ref{WOvSF}, whose proof is more technically delicate --- still hold over this stronger base theory.} We will address these, and other, aspects of the base theory issue in a future paper. However, the search for the ``right" base theory is very fertile mathematical ground, drawing on and responding to foundational ideas from proof theory, generalized recursion theory, and even set theory, and deserves attention from many corners and active debate.

\bibliography{Test}
\bibliographystyle{alpha}

\end{document}